\newif\ifsiam\siamfalse

\PassOptionsToPackage{usenames,dvipsnames,table}{xcolor}
\ifsiam
  \documentclass[review,onefignum,onetabnum]{siamart250211}
\else
  \documentclass[USenglish,10pt,a4paper]{article}
\fi
  \PassOptionsToPackage{noend}{algpseudocode}
  \PassOptionsToPackage{capitalize,nameinlink}{cleveref}
  \PassOptionsToPackage{breakspaces}{clevethm}
  \PassOptionsToPackage{dvipsnames}{xcolor}

  \usepackage[cal=boondoxo, scr=euler]{mathalfa}
  \ifsiam
    \usepackage[
      noload={amsthm},
      excludethm=all,
    ]{myPreamble}
  \else
    \usepackage[bottom=1.18in,left=1.1in,right=1.1in]{geometry}
    \usepackage[
      alglinenumber,
      eqreset=section,
      thmreset=section,
    ]{myPreamble}
  \fi

\makeatletter
  \foreach \x in {gray, blue, orange, red, black} {%
    \edef\temp{%
      \noexpand\expandafter\noexpand\gdef\noexpand\csname\x\noexpand\endcsname{\noexpand\color{\x}}
    }%
    \temp
  }

  \let\OLDand\and
  \def\and{\texorpdfstring{\OLDand}{, }}%
  \ifsiam
    \newcommand{\Sep}{,\ }
  \else
    \newenvironment{keywords}{\par\noindent{\bf Keywords. }}{}
    \newenvironment{AMS}{\par\noindent{\bf AMS subject classifications. }}{}
    \newcommand{\Sep}{ \(\cdot\)\ }
  \fi

  \ifsiam
    \renewenvironment{qedequation*}{\[}{\]}
    \let\qedhere\relax
    
    \theoremstyle{plain}
    \theoremheaderfont{\normalfont\itshape}
    \theorembodyfont{\normalfont}
    \theoremseparator{.}
    \theoremsymbol{}
    \newtheorem{assumption}{Assumption}
    \renewcommand\theassumption{\Roman{assumption}}
    \newtheorem{stepsize}{Stepsize rule}
    
    \newtheorem{relaxation}{Relaxation parameter rule}

    \newsiamremark{example}{Example}
    \newsiamremark{fact}{Fact}
    \newsiamremark{remark}{Remark}

    \CreateTheoremLists{assumption}
    \CreateTheoremLists{corollary}
    \CreateTheoremLists{definition}
    \CreateTheoremLists{example}
    \CreateTheoremLists{fact}
    \CreateTheoremLists{lemma}
    \CreateTheoremLists{proposition}
    \CreateTheoremLists{remark}
    \CreateTheoremLists{theorem}
    \CreateTheoremLists{stepsize}
    \CreateTheoremLists{relaxation}

    \makeatletter
      \let\oldproposition\proposition
      \let\oldendproposition\endproposition
      \renewenvironment{proposition}[1][]{%
        \def\@currentlabelname{#1}\ifstrempty{#1}{\oldproposition}{\oldproposition[#1]}%
      }{\oldendproposition}
      \let\oldlemma\lemma
      \let\oldendlemma\endlemma
      \renewenvironment{lemma}[1][]{%
        \def\@currentlabelname{#1}\ifstrempty{#1}{\oldlemma}{\oldlemma[#1]}%
      }{\oldendlemma}
      \let\oldtheorem\theorem
      \let\oldendtheorem\endtheorem
      \renewenvironment{theorem}[1][]{%
        \def\@currentlabelname{#1}\ifstrempty{#1}{\oldtheorem}{\oldtheorem[#1]}%
      }{\oldendtheorem}
      \let\oldcorollary\corollary
      \let\oldendcorollary\endcorollary
      \renewenvironment{corollary}[1][]{%
        \def\@currentlabelname{#1}\ifstrempty{#1}{\oldcorollary}{\oldcorollary[#1]}%
      }{\oldendcorollary}
    \makeatother

    \foreach \x in {assumption, corollary, definition, example, fact, lemma, proposition, remark, theorem, stepsize, relaxation} {
      \edef\y{\expandafter\euppercase\x}
      \RegisterTheoremName{\x}{\y}
    }
  \else
    \newenvironment{assumption}{\begin{ass}}{\end{ass}}
    \newenvironment{corollary}{\begin{cor}}{\end{cor}}
    \newenvironment{definition}{\begin{defin}}{\end{defin}}
    \newenvironment{example}{\begin{es}}{\end{es}}
    \newenvironment{lemma}{\begin{lem}}{\end{lem}}
    \newenvironment{proposition}{\begin{prop}}{\end{prop}}
    \newenvironment{remark}{\begin{rem}}{\end{rem}}
    \newenvironment{theorem}{\begin{thm}}{\end{thm}}

  \fi

  \allowdisplaybreaks

  \makeatletter
    \renewcommand{\clevethm@proofsectiontitle}{Omitted proofs of }
  \makeatother
  \Crefname{assumption}{Assumption}{Assumptions}
  \Crefname{example}{Example}{Examples}
  \Crefname{fact}{Fact}{Facts}
  \Crefname{remark}{Remark}{Remarks}
  \Crefname{cor}{Corollary}{Corollaries}
  \crefname{ALG@line}{step}{steps}
  \crefname{enumeratpropi}{property}{properties}
  \crefname{enumeratpropii}{property}{properties}
  \counterwithout{algorithm}{section}
  \counterwithout{table}{section}
  \renewcommand{\thealgorithm}{\arabic{algorithm}}
  
  \crefformat{equation}{(#2#1#3)}

  \crefname{table}{table}{tables}
  \Crefname{table}{Table}{Tables}

\usepackage{hyphenat}
\usepackage{cases} 
\usepackage{diagbox}
\usepackage{tabularray}
\UseTblrLibrary{diagbox}
\NewColumnType{M}[1]{Q[m,c,#1]}  
\usepackage{makecell}
\usepackage{makebox}
\usepackage{calrsfs} 
\DeclareMathAlphabet{\pazocal}{OMS}{zplm}{m}{n}

\makeatletter
\newsavebox{\@brx}
\newcommand{\llangle}[1][]{\savebox{\@brx}{\(\m@th{#1\langle}\)}%
  \mathopen{\copy\@brx\mkern2mu\kern-0.9\wd\@brx\usebox{\@brx}}}
\newcommand{\rrangle}[1][]{\savebox{\@brx}{\(\m@th{#1\rangle}\)}%
  \mathclose{\copy\@brx\mkern2mu\kern-0.9\wd\@brx\usebox{\@brx}}}
\makeatother

\makeatletter
\newcommand{\qindef}{\@ifstar{\@qindefb}{\@qindefa}}
\newcommand{\@qindefb}[2]{\inner{#1,#1}_{#2}}
\newcommand{\@qindefa}[2]{\operatorname{q}_{#2}(#1)}
\makeatother

\newcommand{\qsemi}[2]{\operatorname{q}_{#2}(#1)}

\def\mon{\mu}
\def\com{\rho}

\def\Mon{M}
\def\Com{R}

\def\etamin{\bar\eta}

\def\M{P}                           
\def\projR{Q}                       

\def\DRSRho{V}

\def\s1{\|L\|}

\newcommand\optional[1]{[#1]}



\newcommand\tinprod[2]{\langle {#1},\, {#2} \rangle}

\newcommand\other[1]{{#1}^\prime} 

\newcommand{\x}{\bar x}
\newcommand{\y}{\bar y}

\renewcommand\thealgorithm{\arabic{algorithm}}

\makeatletter
  \renewcommand{\appendixproof@toc}{subsubsection}

\makeatother

\newcommand{\dfn}{\coloneqq}
\newcommand{\gph}{\graph}

\DeclarePairedDelimiter{\nrm}{\lVert}{\rVert}
\DeclarePairedDelimiter{\inner}{\langle}{\rangle}

\DeclareMathOperator{\Jac}{J}%

\let\oldBox\Box
\renewcommand{\Box}{\mathbin{\oldBox}}
\newcommand\sym[1]{\mathop{\mathbb{S}}^{#1}}

\pgfplotsset{
  compat=1.17,
  colormap={paper}{rgb255(0cm)=(230,230,230); rgb255(1cm)=(200,200,200); rgb255(2cm)=(185,185,185)}
  }
\usepgfplotslibrary{fillbetween, colormaps}
\usetikzlibrary{arrows.meta, calc, shapes.geometric, perspective}
\usetikzlibrary{decorations.pathreplacing,decorations.markings,arrows,shapes.misc}
\pgfplotscolorbarCMYKworkaroundfalse
\showgrayfalse
\usepackage{caption}

\definecolor{PaperBlue}{HTML}{185477}
\definecolor{PaperGreen}{HTML}{228B22}
\definecolor{PaperOrange}{HTML}{E69C24}
\definecolor{PaperRed}{HTML}{902A3C}

\newif\iftables\tablesfalse
\newif\iftablesCP\tablesCPfalse
\newif\iftablesCPeta\tablesCPetafalse
\tablesCPetatrue

  \newif\ifshowold\showoldfalse
  \newif\ifshownew\shownewfalse
  \AtBeginDocument{%
    \ifshownew\else
      \ifshowold
        \PackageWarning{DioBestia}{WARNING: Showing only OLD version}%
      \else
        \PackageWarning{DioBestia}{ERROR: One at least between \protect\ifshow old and \protect\ifshow new should be True}%
        \PackageError{DioBestia}{One at least between \protect\ifshow old and \protect\ifshow new should be True}{}%
      \fi
    \fi
  }
  
  \pretocmd\grayout{\ifshowgray\else{\gray\hrule\hspace*{0pt}\hfill (removed content)\hfill\hspace*{0pt}\hrule}}{}{}

    \colorlet{newcolor}{orange!70!red}
    \colorlet{oldcolor}{black!30}
    
    \newcommand{\disablecolorlinks}{\def\HyColor@UseColor##1{}}

    \newcommand{\citenum}[1]{{%
      \def\@cite##1##2{{##1\if@tempswa , ##2\fi}}%
      \cite{#1}%
    }}%

  \newcounter{saveTheorem}
  \newcounter{saveEquation}

  \makeatletter
  \newcommand{\mytag}[2]{%
    \text{#1}%
    \@bsphack
    \begingroup
      \@onelevel@sanitize\@currentlabelname
      \edef\@currentlabelname{%
        \expandafter\strip@period\@currentlabelname\relax.\relax\@@@%
      }%
      \protected@write\@auxout{}{%
        \string\newlabel{#2}{%
          {#1}%
          {\thepage}%
          {\@currentlabelname}%
          {\@currentHref}{}%
        }%
      }%
    \endgroup
    \@esphack
  }
  \makeatother
\pgfplotsset{compat=1.16}
\usepgfplotslibrary{fillbetween}
\showgrayfalse

\showoldtrue               
\shownewtrue               
\newcommand{\TheShortTitle}{%
  Spingarn's Method Beyond Elicitable Monotonicity
}
\newcommand{\TheTitle}{
  Spingarn's Method and Progressive Decoupling Beyond Elicitable Monotonicity
}

\newcommand{\TheShortAuthor}{B. Evens, P. Latafat, and P. Patrinos}
\newcommand{\TheFunding}{%
  This work was supported by:
  the Research Foundation Flanders (FWO) PhD grant No. 1183822N and research projects G081222N, G033822N, G0A0920N;
  KU Leuven internal funding: C14/24/103 (Rethinking Transformers Through Duality Principles);
  European Union's Horizon 2020 research and innovation programme under the Marie Skłodowska-Curie grant agreement No. 953348.
  P. Latafat is a member of the Gruppo Nazionale per l'Analisi Matematica, la Probabilit\`a e le loro Applicazioni (GNAMPA - National Group for Mathematical Analysis, Probability and their Applications) of the Istituto Nazionale di Alta Matematica (INdAM - National Institute of Higher Mathematics).%
}  
\newcommand{\TheKeywords}{%
  Spingarn's method of partial inverses\Sep
  progressive decoupling\Sep
  linkage problems\Sep
  nonconvex optimization\Sep
  nonmonotone variational inequalities
}
\newcommand{\TheAMSsubj}{%
  47H04\Sep 
  49J52\Sep 
  49J53\Sep
  65K15\Sep 
  90C26. 
}
\newcommand{\TheAbstract}{%
  Spingarn's method of partial inverses and the progressive decoupling algorithm address inclusion problems involving the sum of an operator and the normal cone of a linear subspace, known as linkage problems.
  Despite their success, existing convergence results are limited to the so-called elicitable monotone setting, where nonmonotonicity is allowed only on the orthogonal complement of the linkage subspace.
  In this paper, we introduce progressive decoupling+, a generalized version of standard progressive decoupling that incorporates separate relaxation parameters for the linkage subspace and its orthogonal complement.
  We prove convergence under conditions that link the relaxation parameters to the nonmonotonicity of their respective subspaces and show that the special cases of Spingarn's method and standard progressive decoupling also extend beyond the elicitable monotone setting.
  Our analysis hinges upon an equivalence between progressive decoupling+ and the preconditioned proximal point algorithm, for which we develop a general local convergence analysis in a certain nonmonotone setting.
}

\author{%
    Brecht Evens\textsuperscript{1}
    \and
    Puya Latafat\textsuperscript{2}
    \and
    Panagiotis Patrinos\textsuperscript{1}
    \thanks{
      \TheFunding\newline
      \indent
      \textsuperscript{1}\TheAddressKU,
      {\tt
        \{%
          \href{mailto:brecht.evens@kuleuven.be}{brecht.evens},%
          \href{mailto:panos.patrinos@kuleuven.be}{panos.patrinos}%
          \}\texttt{@kuleuven.be}.
      }\newline
      \indent
      \textsuperscript{2}\TheAddressIMT, 
      {\tt
          \href{mailto:puya.latafat@imtlucca.it}{puya.latafat@imtlucca.it}.
      }
    }
  }
\renewcommand\footnotemark{}

\ifsiam
  \headers{\TheShortTitle}{\TheShortAuthor}
  
  \title{%
  \TheTitle%
  \thanks{%
    \TheAddressKU\newline
    \textit{E-mail:}
    {\tt
      \{%
        \href{mailto:brecht.evens@kuleuven.be}{brecht.evens},%
        \href{mailto:puya.latafat@imtlucca.it}{puya.latafat},%
        \href{mailto:panos.patrinos@kuleuven.be}{panos.patrinos}%
        \}@kuleuven.be.
    }%
    Submitted to the editors \today.%
    \funding{\TheFunding}%
  }%
}

\else
  \title{\texorpdfstring{\TheTitle}{\TheShortTitle}}
  \date{}
\fi

\newif\ifCPlong\CPlongfalse
\CPlongtrue
\newif\ifshowfullCP\showfullCPtrue

\begin{document}
    \maketitle

    \begin{abstract}
      \TheAbstract
    \end{abstract}
    

\section{Background}

A fundamental problem in variational analysis and optimization is to solve structured inclusion problems of the form
\begin{equation}\label{prob:inclusion-intro}
	\text{find} \quad x \in \R^n \quad \text{such that} \quad 0 \in N_X(x)+ S(x),
\end{equation}
where
\(S : \R^n \rightrightarrows \R^n\)
is a set-valued operator, 
\(X \subset \R^n\) is a proper closed linear subspace and \(N_X\) is the normal cone operator of $X$.
Due to the linearity of the subspace $X$,
this problem can be equivalently represented as a so-called \emph{linkage problem} \cite{rockafellar2019Progressive}, given by
\begin{equation}\label{prob:linkage-intro}
	\text{find} \quad (x, y) \in \R^n \times \R^n \quad \text{such that} \quad x \in X, \; y\in X^{\bot},\; y\in S(x).\tag{P}
\end{equation}
The condition $x \in X$ is known as the \emph{linkage constraint} with corresponding multipliers $y \in X^\bot$. 
Throughout, we denote to the set of all pairs \((x^\star, y^\star)\) solving \eqref{prob:linkage-intro} as \( \link_X S \).

Linkage problems arise naturally in many applications as first-order conditions of optimality. 
As an important special case consider the following class of constrained optimization problems, known as \emph{monotropic programming} problems  
\cite{rockafellar1999network,bertsekas2008extended}:
\begin{equation}\label{prob:linkage-optimization}
	\minimize_{x\in\R^{n}} \quad \psi(x) \quad \stt \quad x \in X,
\end{equation}
where $\psi : \R^n \rightarrow \Rinf \coloneqq \R\cup\set\infty$ is a proper lower semicontinuous (lsc) function.
Under suitable constraint qualifications on the function $\psi$, the first-order optimality conditions of \eqref{prob:linkage-optimization} take the form of linkage problem \eqref{prob:linkage-intro} with $S$ being to the subdifferential $\partial \psi$.

Many optimization problems can be cast as an instance of \eqref{prob:linkage-optimization} by introducing auxiliary variables and using the subspace $X$ to enforce relationships between these optimization variables.
For instance, consider the structured minimization problem 
\[
	\minimize_{x\in\R^{n}}
	\quad
	g(x) + h(Lx)
\]
where \(g:\R^n\to \Rinf\) and \(h:\R^m\to \Rinf\) are lsc functions and \(L:\R^n \to \R^m\) is a linear mapping. Such problems can be equivalently formulated as \eqref{prob:linkage-optimization} using a simple lifting, defining the function \(\psi(x_1,x_2) = g(x_1)+h(x_2)\) and linear subspace \(X = \set{(v_1,v_2)\in \R^{n+m}}[v_2 = Lv_1]\). 
Other examples of optimization problems of the form \eqref{prob:linkage-optimization} arise in Lagrangian duality \cite[\S 29, \S30]{Rockafellar1970Convex}, \cite[\S 11.H]{rockafellar2009Variational} and in multistage stochastic programming \cite{ruszczynski2003stochastic,de2021risk}, where the subspace \(X\) enforces non-anticipativity in the decision-making.

However, the range of applications is not only restricted to minimization problems.
For instance, consider the operator splitting problem of finding a vector \(x\in \R^n\) such that
\[
	0 \in A_1(x) + A_2(x) +\ldots + A_N(x),
\]
where \(A_i : \R^n \rightrightarrows \R^n\) are general set-valued operators. Such structures capture the popular finite sum minimization problem as well as matrix splitting problems. It can be cast in the form of linkage problem \eqref{prob:linkage-intro} by letting
\(
	S(x_1,\ldots,x_N)
		=
	(A_1(x_1),\ldots,A_N(x_N))
\)
and
\(
	X \subset \R^{Nn}
\)
denoting the consensus set (see e.g. \Cref{ex:consensus:cond,ex:linear-system}). 
For more details and examples of linkage problems, we refer the interested reader to \cite[\S 4]{rockafellar2019Progressive}.

The most popular method for solving \eqref{prob:linkage-intro} is Spingarn's method of partial inverses, which involves a single evaluation of the resolvent $J_S \coloneqq (\id + S)^{-1}$ at each iteration, followed by projections onto the subspaces $X$ and $X^\bot$, respectively.
In particular, given an initial guess \(x^0 \in X\) and \(y^0 \in X^\bot\), the iterations are given by
\begin{equation}
	\begin{cases}\tag{Spingarn}\label{eq:spingarn}
		\displaystyle   
		q^k &\in J_{S}(x^k + y^k)\\
		x^{k+1}&= \proj_X (q^k)\\
		y^{k+1}&= y^{k} - \proj_{X^\bot} (q^k),
	\end{cases}
\end{equation}
where $\proj_X$ and $\proj_{X^\bot}$ denote projections onto $X$ and $X^\bot$, respectively.
Under the standard assumption that $S$ is maximally monotone, Spingarn showed that the sequence $\seq{x^k, y^k}$ generated by this method converges to some point in \(\link_X S\), provided one exists \cite{spingarn1983Partial}.

More recently, Rockafellar proposed a generalization of Spingarn's method of partial inverses, named the \emph{progressive decoupling} algorithm, which introduces a stepsize parameter $\gamma$ and a so-called \emph{elicitation} parameter $\mon \leq 0$ into the update rule as follows \cite{rockafellar2019Progressive}:
\begin{equation}
	\begin{cases}
		\displaystyle   
		q^k &\in J_{\gamma^{-1} S}(x^k + \gamma^{-1}y^k)\\
		x^{k+1}&= \proj_X (q^k)\\
		y^{k+1}&= y^{k} - (\gamma + \mon)\proj_{X^\bot} (q^k).
	\end{cases}
	\label{eq:progdec-standard}
\end{equation}
This algorithm can be viewed as a generalization of the \emph{progressive hedging algorithm} for multistage stochastic programs \cite{rockafellar1991scenarios} to general linkage problems.
The progressive decoupling algorithm reduces to Spingarn's method of partial inverses for $\mon = 0$ and $\gamma = 1$.
By drawing a connection with the proximal point method, Rockafellar established that progressive decoupling converges under the assumption that $S$ is \emph{elicitable monotone}, i.e., that
\[
	\text{
		\(S-\mon\proj_{X^\bot}\) is monotone for some
		\(\mon \leq 0\).
		}
\]
In this work, we introduce a further generalization of Spingarn's method, hereafter referred to as \emph{progressive decoupling+}, for solving linkage problems \eqref{prob:linkage-intro} in the nonmonotone setting.
It generalizes the standard progressive decoupling algorithm by incorporating two relaxation parameters, $\lambda_x$ and $\lambda_y$, whose main role is to allow for greater flexibility in handling the nonmonotonicity of the operator $S$.
Specifically, given an initial guess \(x^0 \in X\) and \(y^0 \in X^\bot\), the iterates are as follows:
\begin{equation}
	\begin{cases}
		\displaystyle   
		q^k &\in J_{\gamma^{-1} S}(x^k + \gamma^{-1}y^k)\\
		x^{k+1}&= (1 - \lambda_{x})x^{k} + \lambda_{x}\proj_X (q^k)\\
		y^{k+1}&= y^{k} - \lambda_{y}\gamma\proj_{X^\bot} (q^k).
	\end{cases}
	\tag{ProgDec+}\label{eq:progdec}
\end{equation}
Leveraging the full generality of this algorithm, we establish the convergence of \ref{eq:progdec} under the assumption that
\begin{align*}
	\text{
		\(
			(S - \mon \proj_{X^\bot})^{-1} - \com \proj_X	
		\)
		is monotone for some
		\(\mon \in \R\)
		and
		\(\com \in \R\),
		}
	\numberthis\label{eq:intro:ass-progdec}
\end{align*}
a condition referred to as \emph{\((\mon \proj_{X^\bot}, \com \proj_X)\)\hyp{}semimonotonicity} of \(S\) (see \cref{def:semimonotonicity}).
It is not required that this monotonicity assumption holds between any two points in the graph of $S$, but instead only between any point in the graph of $S$ and at least one element of the solution set $\link_X S$ (see \cref{ass:progdec-local}).
This operator class arises quite naturally in practical applications, as we illustrate with several examples throughout the paper (see for instance \Cref{ex:Rosenbrock:progdec,ex:double-well,ex:consensus:cond} and the calculus rules for continuously differentiable mappings in \Cref{subsec:progdec:examples}).
Notably, it strictly encompasses both the elicitable monotone setting (retrieved for \(\com = 0\)) and the standard monotone setting (retrieved for \(\mon = \com = 0\)). For other recent works studying the class of semimonotone operators, we refer to \cite{evens2023convergence,evens2023convergenceCP,quan2024scaled}.

\ref{eq:progdec} recovers several existing methods as special cases. Notably, when $\lambda_x = \lambda_y$, it reduces to the well-known relaxed Douglas--Rachford splitting (DRS) method applied to inclusion problem \eqref{prob:inclusion-intro} (see \Cref{rem:equiv-progdec-DRS}).
The standard progressive decoupling method is recovered for \(\mon \leq 0\), \(\lambda_x = 1\), and \(\lambda_y = 1 + \nicefrac{\mon}{\gamma}\), while Spingarn's method corresponds to the choice \(\gamma = \lambda_x = \lambda_y = 1\). 
Our novel convergence theory extends existing results for these three well-known methods, as summarized in \Cref{tab:progdec:connections}, reaching beyond the traditional monotone setting.
Specifically, our analysis reveals that the range of admissible relaxation parameters $\lambda_x$ is determined by the properties of $S$ on $X$, quantified by the parameter $\com$, while the range of viable $\lambda_y$ is similarly governed by $\mon$.
In particular, higher levels of nonmonotonicity correspond to a narrower range of permissible relaxation parameters (see equation \eqref{eq:thm:progdec-local:stepsize}).

Our development of \ref{eq:progdec} primarily relies on two key equivalences, with the \emph{partial inverse} of the operator $S$ playing a central role (see \cref{def:partial}).
Specifically, we first show that \ref{eq:progdec} can be interpreted as an instance of the preconditioned proximal point algorithm (PPPA) with matrix  
\(
	\lambda_x \proj_X + \lambda_y \proj_{X^\bot}
\)
as its relaxation parameter applied to the partial inverse of $S$. Second, we demonstrate that under the semimonotonicity assumption on $S$ from \eqref{eq:intro:ass-progdec}, the partial inverse of $S$ satisfies an \emph{oblique weak Minty assumption} \cite{evens2023convergence} (see \Cref{def:SWMVI}).
In order to leverage these equivalences, we develop a general convergence analysis for relaxed PPPA in the nonmonotone setting with arbitrary relaxation matrices. 
This is achieved by exploiting the technique developed in \cite{evens2023convergence} for nonmonotone inclusion problems that relies on interpreting each iteration of the algorithm as a projection onto a certain hyperplane. 
The hyperplane projection interpretation of the proximal point method is a powerful idea that dates back to \cite{solodov1996modified,solodov1999hybrid,konnov1997Class} in the monotone setting.

To further broaden the scope of our work, our convergence analysis includes the local setting, where the assumptions only hold on a subset of the corresponding graph.
Specifically, we provide convergence results for relaxed PPPA applied to inclusion problems that admit \emph{local} oblique weak Minty solutions (see \cref{def:SWMVI}).
The idea of graph localization draws on Pennanen's local analysis of the proximal point algorithm \cite{pennanen2002Local,iusem2003Inexact,combettes2004proximal}, where a key step in the analysis is to ensure that the PPPA iterates remain within the localized, restricted graph.
In the minimization setting, this localization of the graph is closely related to relaxing global convexity to \emph{variational convexity} \cite{rockafellar2019varconv,rockafellar2023augmented}.
Note that many nonconvex optimization problems possess such a variational convexity assumption, providing a partial explanation for the success of many numerical methods beyond the standard, monotone case.
We would like to highlight that our analysis is valuable not only for studying the convergence behavior of \ref{eq:progdec}, but also in its own right, as many other methods are known to be instances of relaxed PPPA, including the alternating direction method of multipliers (ADMM) \cite[\S 8]{eckstein1988LionsMercier}, the (proximal) augmented Lagrangian method (ALM) \cite[\S 4 \& 5]{rockafellar1976augmented} and the primal-dual hybrid gradient (PDHG) method \cite[Lem. 2.2]{he2012Convergence}, \cite{esser2010general}, also known as Chambolle-Pock (CP) \cite{chambolle2011firstorder}.


  \subsection{Organization}

The paper is structured as follows. First, we recall some notation in \Cref{subsec:notation}.
Then, we demonstrate in \Cref{sec:equivalence-pppa-progdec} that the proposed algorithm can be cast as a specific instance of the relaxed preconditioned proximal point algorithm.
Consequently, we dedicate \Cref{sec:local-pppa} to the local convergence of relaxed PPPA under an oblique weak Minty condition.
\Cref{sec:progdec} then applies these convergence results to our algorithm, showing local convergence to a solution of linkage problem \eqref{prob:linkage-intro} under the assumption that $S$ is locally semimonotone.
Moreover, several examples are provided, including an example showing tightness of the obtained range for the relaxation parameters, and we provide sufficient conditions for a differentiable operator to be semimonotone.
Finally, \cref{sec:conclusion} concludes the paper.

  \subsection{Notation}\label{subsec:notation}

The set of natural numbers including zero is denoted by \(\N\coloneqq\set{0,1,\hdots}\).
The set of real and extended-real numbers are denoted by \(\R\coloneqq(-\infty,\infty)\) and \(\Rinf\coloneqq\R\cup\set\infty\), while the positive and strictly positive reals are \(\R_+\coloneqq[0,\infty)\) and \(\R_{++}\coloneqq(0,\infty)\).
We use the notation $\seq{w^k}[k\in I]$ to denote a sequence with indices in the set $I\subseteq \N$. When dealing with scalar sequences we use the subscript notation $\seq{\gamma_k}[k\in I]$. 
We denote the positive part of a real number by $[\cdot]_{+} \dfn \max\{0, \cdot\}$ and the negative part by $[\cdot]_{-} \dfn -\min\{0, \cdot\}$.
With \(\id\) we indicate the identity function \(x\mapsto x\) defined on a suitable space.
The identity matrix is denoted by $\I_n\in\R^{n\times n}$ and the zero matrix by $0_{m\times n}\in\R^{m\times n}$; we write respectively $\I$ and $0$ when no ambiguity occurs. 

We denote by $\R^n$ the standard $n$-dimensional Euclidean space with inner product $\langle\cdot,\cdot\rangle$ and induced norm $\|\cdot\|$. 
The set of symmetric $n$-by-$n$ matrices is denoted by $\sym{n}$. Given a symmetric matrix $P \in \sym{n}$, we write $P \succeq 0$ and $P \succ 0$ to denote that $P$ is positive semidefinite and positive definite, respectively.
Furthermore, for any $P \in \sym{n}$ we define the 
quadratic function \(\qsemi{x}{P} \dfn \inner{x, Px}\).
Let $\diag(\cdot)$ denote the diagonal
matrix whose arguments constitute its diagonal elements.
For arbitrary matrices $A$ and $B$, we define the direct sum $A \oplus B = \blkdiag(A,B)$,
where
$\blkdiag(\cdot)$ denotes the block diagonal
matrix whose arguments constitute its diagonal blocks.

Two vectors $u, v \in \R^n$ are said to be orthogonal if $\inner{u,v} = 0$, and orthonormal if they are orthogonal and $\nrm{u} = \nrm{v} = 1$. Two linear subspaces $\mathbf{U} \subseteq \R^n$ and $\mathbf{V} \subseteq \R^n$ are said to be orthogonal if any $u \in \mathbf{U}$ and any $v \in \mathbf{V}$ are orthogonal.
We say that $U \in \R^{n \times m}$ is an orthonormal basis for a linear subspace $\mathbf{U} \subseteq \R^n$ if $U$ has orthonormal columns and 
the range of $U$ is given by $\mathbf{U}$.
 

The domain of an extended real-valued function $f : \R^n \rightarrow \Rinf$ is the set $\dom f \dfn \set{x \in \R^n}[f(x) < \infty]$.
We say that $f$ is proper if $\dom f \neq \emptyset$ and that $f$ is lower semicontinuous (lsc) if the epigraph $\epi f \dfn \set{(x,\alpha) \in \R^n \times \R}[f(x) \leq \alpha]$ is a closed subset of $\R^{n+1}$.
We denote the limiting subdifferential of $f$ by $\partial f$.
We denote the normal cone of a set $E\subseteq\R^n$ by \(N_E\).
The projection onto and the distance from $E$ with respect to $\|\cdot\|_{\projR}$, \(Q\in\sym{n}_{++}\), are denoted by 
\(
    \proj^Q_E(x)\coloneqq{}\argmin_{z\in E}\{\|z-x\|_{\projR}\}
\)
and
\(
    \dist_{\projR}(x,E)\coloneqq{} \inf_{z\in E}\{\|z-x\|_{\projR}\},
\)
respectively.
For two sets \( E, F \subseteq \mathbb{R}^n \), we denote the Minkowski sum by
\(
    E + F \coloneqq \{e + f \mid e \in E, \, f \in F\}
\).
An operator or set-valued mapping $A:\R^n\rightrightarrows\R^d$ maps each point $x\in\R^n$ to a subset $A(x)$ of $\R^d$. We will use the notation $A(x)$ and $Ax$ interchangeably. 
We denote the domain of $A$ by $\dom A\coloneqq\{x\in\R^n\mid Ax\neq\emptyset\}$,
its graph by $\graph A\coloneqq\{(x,y)\in\R^n\times \R^d\mid y\in Ax\}$, and 
the set of its zeros by $\zer A\coloneqq\{x\in\R^n \mid 0\in Ax\}$. 
The inverse of $A$ is defined through its graph: $\graph A^{-1}\coloneqq\{(y,x)\mid (x,y)\in\graph A\}$.
The \emph{resolvent} of $A$ is defined by $J_A\coloneqq(\id+A)^{-1}$.
The composition of two operators \(A\) and \(B\) is denoted by \(A \circ B\).
We say that \(A\) is \emph{outer semicontinuous (osc)} at $\other{x}\in\dom A$ if 
\ifsiam
\else
    \begin{equation*}
        \limsup_{x\to \other{x}} Ax \coloneqq \{y\mid \exists x^k \to \other{x}, \exists y^k \to y \textrm{ with } y^k\in Ax^k\} \subseteq A\other{x}.
    \end{equation*}
\fi
\ifsiam
\else
    Outer semicontinuity of $A$ everywhere is equivalent to its graph being a closed subset of $\R^n\times \R^d$.
\fi

\begin{definition}[local $\DRSRho$\hyp{}comonotonicity]\label{def:comon}
    Let $\DRSRho\in\R^{n \times n}$ be a symmetric (possible indefinite) matrix.
    An operator \(T : \R^n \rightrightarrows \R^n\) is said to be $\DRSRho$\hyp{}comonotone at \((x', y') \in \gph T\) 
    on a set $\pazocal U \ni (x', y')$
    if 
    \begin{equation*}
        \label{eq:def:comon}
        \langle x-x',y-y'\rangle\geq \qindef{y-y'}{\DRSRho},\qquad \text{for all $(x, y) \in \gph T \cap \pazocal U$},
    \end{equation*}
    where the quadratic form \(\qindef{\cdot}{\DRSRho} \dfn \langle \cdot, \cdot \rangle_{\DRSRho}\).
    An operator \(T\) is said to be $\DRSRho$\hyp{}comonotone on $\pazocal U$ if it is $\DRSRho$\hyp{}comonotone at all $(\other{x}, \other{y})\in \gph T \cap \pazocal U$.
    It is said to be maximally $\DRSRho$\hyp{}comonotone on $\pazocal U$ if its graph is not strictly contained in the graph of another $\DRSRho$\hyp{}comonotone operator on $\pazocal U$. 

    Throughout, whenever $\pazocal U = \R^n \times \R^n$, the set $\pazocal U$ is omitted, and whenever $\DRSRho = \rho \I$ for some $\rho \in \R$, the prefix $\DRSRho$ is replaced by $\rho$ and condition \eqref{eq:def:comon} reduces to
    \begin{equation*}
        \langle x-x',y-y'\rangle\geq \rho\|y-y'\|^2,\qquad \text{for all $(x, y) \in \gph T \cap \pazocal U$},
    \end{equation*}
\end{definition}

\section{Progressive decoupling+ as an instance of relaxed PPPA}\label{sec:equivalence-pppa-progdec}

We recall the definition of the partial inverse of an operator, introduced in \cite{spingarn1983Partial}, along with some calculus rules which will be used throughout the paper. The partial inverse will play a central role in our upcoming analysis. For more details on the partial inverse, we refer the interested reader to \cite[\S 20.3]{bauschke2017Convex}.

\begin{definition}\label{def:partial}
	Let \(X\subseteq \R^n\) be a closed linear subspace, and \(T:\R^n\rightrightarrows\R^n\) be a set-valued operator. Then, the partial inverse of \(T\) with respect to \(X\) is the operator \(T^X:\R^n\rightrightarrows\R^n\) defined by 
	\begin{equation}\label{eq:def:partial}
		\graph T^X \coloneqq L_X(\graph(T)),
	\end{equation}
	where \(L_X\)  is the Spingarn operator
	\(
		L_X (x,y) 
			{}\coloneqq{}
		\big(
		\proj_X(x) + \proj_{X^\bot}(y),
		\proj_X(y) + \proj_{X^\bot}(x)
		\big). 
	\)
\end{definition}

\begin{fact}
    Let \(X\subseteq \R^n\) be a closed linear subspace, and \(T:\R^n\rightrightarrows\R^n\) be a set-valued operator.
    Then, the following hold.
    \begin{enumerate}
        \item\label{it:spingarn:inverse} \(T^X = (T^{X^\bot})^{-1} = (T^{-1})^{X^\bot}\).
        \item\label{it:spingarn:gph} \(L_X^{-1} = L_X\) and \(\graph T = L_X(\graph(T^X))\).
        \item\label{it:spingarn:zero} \(z \in \zer T^X\) if and only if \((\proj_X(z), \proj_{X^\bot}(z)) \in \link_X T\).
		\item\label{it:spingarn:linear} If \(T \in \R^{n \times n}\) and $\proj_X + \proj_{X^\bot} T$ is invertible, then
		\(
			T^X
				{}={}
			(\proj_{X^\bot} + \proj_X T)(\proj_X + \proj_{X^\bot} T)^{-1}
		\). 
    \end{enumerate}
\end{fact}
Notably, \Cref{it:spingarn:zero} demonstrates that 
a pair $(x, y) \in \link_X S$ can be obtained by finding a zero of the operator $S^X$.
Our upcoming analysis will revolve around this particular equivalence.
Specifically, we introduce the 
inclusion problem of finding a zero of a set-valued operator \(T:\R^n \rightrightarrows \R^n\), \ie,
\begin{equation} \label{prob:P1} \tag{G-I} 
    \text{find} \quad z \in \R^n \quad \text{such that} \quad 0 \in T z. 
\end{equation}
The most well-known algorithm for solving this class of problems is the proximal point algorithm and its preconditioned/relaxed variants. 
Given a symmetric positive definite preconditioning matrix \(\M \in \R^{n\times n}\) and a relaxation matrix \(\Lambda \in \R^{n\times n}\), the (relaxed) preconditioned proximal point algorithm applied to \eqref{prob:P1} consists of the following fixed point iterations:
\begin{equation}\label{eq:PPPA-intro}\tag{PPPA}
    \begin{cases}
        \bar{z}^k&\in(\M+T)^{-1}\M z^k\\
        z^{k+1}&=z^k+\Lambda(\bar{z}^k-z^k).
    \end{cases}
\end{equation}
In what follows, we will show that \ref{eq:progdec} can be cast as a particular instance of the (relaxed) preconditioned proximal point algorithm, as a consequence of the following lemma.

\begin{lemma}\label{lem:progdec:resolvent-equivalence}
	Consider an operator \(S:\R^n\rightrightarrows\R^n\) and a closed linear subspace \(X\subseteq \R^n\).
	Let $\gamma \in (0, +\infty)$ and define \(\M = \gamma \proj_X + \gamma^{-1} \proj_{X^\bot}\).
	Then, for any $z \in \R^n$ it holds that
	\(
		\bar z \in (\M + S^X)^{-1} \M z
	\)
	if and only if
	\[
		\Bigl(\underbrace{\proj_{X}(\bar{z})+\gamma^{-1}\proj_{X^{\bot}}(z-\bar{z})}_{\eqqcolon\, q}, \proj_{X^{\bot}}(\bar{z}) + \gamma\proj_{X}(z-\bar{z})\Bigr)
			{}\in{} 
		\gph S,
	\]
	or equivalently,
	\[
		q 
			\coloneqq 
		\proj_{X}(\bar{z}) 
			+
		\gamma^{-1}\proj_{X^{\bot}}(z-\bar{z}) 
			\in 
		J_{\gamma^{-1}S}\bigl(\proj_X(z) + \gamma^{-1}\proj_{X^\bot}(z)\bigr).
	\]
	\begin{proof}
		If 
		\(
			\bar z \in (\M + S^X)^{-1} \M z
		\)
		then it holds that \((\bar z, \M (z - \bar z)) \in \graph S^X\), so that equivalently \(L_X(\bar z, \M (z - \bar z)) \in L_X(\graph S^X) = \graph S\) owing to \cref{it:spingarn:gph}. In particular, 
		\[
		\big(
			\proj_{X}(\bar{z})+\proj_{X^{\bot}}(\M (z-\bar{z})),
			\proj_{X^{\bot}}(\bar{z})+\proj_{X}(\M (z-\bar{z}))
		\big)
			{}\in{} 
		\graph S.
		\]
		Since \(\proj_X \circ \M = \gamma \proj_X\) and \(\proj_{X^\bot} \circ \M = \gamma^{-1} \proj_{X^\bot}\), this corresponds to
		\begin{align*}
			\proj_{X^{\bot}}(\bar{z}) + \gamma\proj_{X}(z-\bar{z})
				{}\in{} 
			S(\underbrace{\proj_{X}(\bar{z})+\gamma^{-1}\proj_{X^{\bot}}(z-\bar{z})}_{\eqqcolon\, q}).
		\end{align*}
		Multiplying by $\gamma^{-1} > 0$ and adding $q$ to both sides, we obtain that
		\begin{equation}\label{eq::progdec:resolvent-equivalence}
			\proj_{X}(z) + \gamma^{-1}\proj_{X^{\bot}}(z)
				{}\in{} 
			(\id + \gamma^{-1}S)(q).
		\end{equation}
		Since $J_{\gamma^{-1}S} \coloneqq (\id + \gamma^{-1}S)^{-1}$, the proof is completed.
	\end{proof}
\end{lemma}
In light of the above lemma, by selecting the preconditioner
\(
	\M = \gamma \proj_X + \gamma^{-1} \proj_{X^\bot}
\),
relaxation matrix
\(
	\Lambda \coloneqq \lambda_x \proj_X + \lambda_y \proj_{X^\bot}
\)
and denoting $x^k = \proj_X(z^k)$, $y^k = \proj_{X^\bot}(z^k)$, $\bar{x}^k = \proj_X(\bar{z}^k)$ and $\bar{y}^k = \proj_{X^\bot}(\bar{z}^k)$, update rule \eqref{eq:PPPA-intro} reduces to
\begin{equation}
	\label{eq:progdec-full}
	\begin{cases}
		\left.
		\begin{aligned}
			&\makebox[0pt][l]{$q^k$}\phantom{x^{k+1}} 
			\makebox[5.25cm][l]{$
					{}\in{}
				J_{\gamma^{-1} S}(x^k + \gamma^{-1} y^k)
			$}\\
			&\makebox[0pt][l]{$\bar{x}^{k}$}\phantom{x^{k+1}} 
				{}={}
			\proj_{X}(q^k),\\
			&\makebox[0pt][l]{$\bar{y}^{k}$}\phantom{x^{k+1}} 
				{}={}
			y^k - \gamma\proj_{X^{\bot}}(q^k)
		\end{aligned} \color{gray}\right\} & \quad \textcolor{gray}{\text{proximal step } \bar{z}^{k} \in (P + S^X)^{-1} P z^k}\\
		\left.
		\begin{aligned}
			&x^{k+1}
			\makebox[5.25cm][l]{$
					{}={}
				x^k + \lambda_x(\bar{x}^k - x^k)
			$}\\
			&y^{k+1}
				{}={}
			y^k + \lambda_y(\bar{y}^k - y^k)
		\end{aligned} \color{gray}\right\} & \quad \textcolor{gray}{\text{relaxation step } z^{k+1} = z^k + \Lambda(\bar{z}^k - z^k)}
	\end{cases}
\end{equation}
which is exactly \ref{eq:progdec}. 
This equivalence is stated more formally below.
\begin{lemma}[equivalence of \ref{eq:progdec} and \ref{eq:PPPA-intro}]\label{lem:equivalence:PPPA-progdec}
	Let \(x^0 \in X\), \(y^0 \in X^\bot\) and set \(z^0 = x^0 + y^0\).
	Then, to any sequence $(q^k, x^k, y^k)_{k\in\N}$ generated by \ref{eq:progdec} (initialized with $x^0$ and $y^0$) for solving \eqref{prob:linkage-intro} with stepsizes $\gamma, \lambda_x, \lambda_y > 0$ there corresponds a sequence 
	\((z^k, \bar z^k)_{k\in\N} = (x^k + y^k, \proj_X(q^k) + y^k - \gamma\proj_{X^\bot}(q^k))_{k\in\N}\)
	generated by 
	\ref{eq:PPPA-intro} (initialized with $z^0 = x^0 + y^0$) applied to $0 \in S^X(z)$ with 
	\(\M = \gamma \proj_X + \gamma^{-1} \proj_{X^\bot}\)
	and
	\(\Lambda = \lambda_x \proj_X + \lambda_y \proj_{X^\bot}\).
\end{lemma}

In the upcoming section, we will first examine the (local) convergence properties of \ref{eq:PPPA-intro} in the nonmonotone setting. 
Then, in \cref{sec:progdec}, we will 
leverage the equivalence from \Cref{lem:equivalence:PPPA-progdec}
to obtain local convergence results for \ref{eq:progdec} for nonmonotone operators $S$.

\begin{remark}[relaxed Douglas--Rachford splitting]\label{rem:equiv-progdec-DRS}
	We remark that the relaxed variant of the Douglas--Rachford splitting (DRS) method can be seen as a particular instance of our algorithm.  
	Specifically, for a given positive stepsize \( \gamma \) and relaxation parameter \( \lambda \),  
	the relaxed DRS method applied to \eqref{prob:inclusion-intro} follows the iteration:
	\begin{equation*}
		\begin{cases}
			\displaystyle   
			x^k &= \proj_X(s^k)\\
			q^k &\in J_{\gamma^{-1} S}(2x^k - s^k)\\
			s^{k+1} &= s^k + \lambda(q^k - x^k).
		\end{cases}
	\end{equation*}
	To establish the equivalence with \ref{eq:progdec}, introduce the variable  
	\(
		y^k = \gamma(x^k - s^k),
	\)
	and reorder the updates by starting with the update for \( q^k \).  
	This leads to the equivalent formulation:
	\begin{equation*}
		\begin{cases}
			\displaystyle   
			q^k &\in J_{\gamma^{-1} S}(x^k + \gamma^{-1} y^k)\\
			x^{k+1} &= 
			\proj_X(x^k - \gamma^{-1} y^k + \lambda(q^k - x^k))\\
			y^{k+1} &= 
			y^k + \gamma(x^{k+1} - (1 - \lambda) x^k - \lambda_k q^k).
		\end{cases}
	\end{equation*}
	Therefore, if the algorithm is initialized with \( x^0 \in X \) and \( y^0 \in X^\bot \), then it is an instance of \ref{eq:progdec} with \( \lambda_x = \lambda_y \).
	This equivalence extends those from \cite[\S 5]{eckstein1992DouglasRachford} for Spingarn's method (where \(\gamma = \lambda_x = \lambda_y = 1\)) and from \cite[Appendix]{de2021risk} for standard progressive decoupling (where $\mon = 1$ and \(\lambda_x = \lambda_y = 1\)).
\end{remark}

\section{Local convergence of nonmonotone PPPA}\label{sec:local-pppa}

Our analysis in this section builds upon and extends the one developed recently in 
\cite{evens2023convergence} for nonmonotone problems which admit so-called oblique weak Minty solutions. 
We generalize the definition of weak Minty solutions to the local setting, only requiring the involved inequality to hold on a subset of the graph and allow for matrix relaxation parameters, a generalization that is crucial for developing \ref{eq:progdec} as established in \Cref{lem:equivalence:PPPA-progdec}. 
This localization of the graph is inspired by the local analysis for the proximal point algorithm in the maximally comonotone setting \cite{pennanen2002Local,iusem2003Inexact,combettes2004proximal} and related work on variational convexity \cite{rockafellar2019varconv}.

\begin{definition}[local $\DRSRho$\hyp{}oblique weak Minty solutions]\label{def:SWMVI}
    Let $\DRSRho\in\R^{n \times n}$ be a symmetric (possibly indefinite) matrix.
    An operator \(T : \R^n \rightrightarrows \R^n\) is said to have $\DRSRho$\hyp{}oblique weak Minty solutions at a nonempty, convex set \(\pazocal{S}^\star\subseteq\zer T\) 
    on a set $\pazocal U \subseteq \R^n \times \R^n$
    if 
    \begin{equation}\label{eq:def:WMVI}
        \langle v,\bar z-z^\star\rangle\geq \qindef{v}{\DRSRho},\qquad \text{for all $z^\star\in\pazocal{S}^\star$ and all $(\bar z,v) \in \gph T \cap \pazocal U$}, 
    \end{equation}
    where the quadratic form \(\qindef{v}{\DRSRho} \dfn \langle v, \DRSRho v \rangle\).
    Whenever $\DRSRho = \rho \I$ for some $\rho \in \R$, we refer to them as $\rho$\hyp{}weak Minty solutions.
    In the global setting, i.e., when $\pazocal U = \R^n \times \R^n$, the set $\pazocal U$ is omitted.
\end{definition}
It is important to emphasize that even in the global setting, \ie, when $\pazocal U = \R^n \times \R^n$, this definition covers many classical assumptions employed in literature for both proximal and gradient-type methods.
For instance, when $\pazocal U = \R^n \times \R^n$ and $\DRSRho$ is equal to the zero matrix, \eqref{eq:def:WMVI} reduces to the classic Minty variational inequality (MVI) \cite{minty1962, giannessi1998Minty}.
This condition is sometimes also referred to as variational coherence and is satisfied for instance by all pseudoconvex, star-convex and quasar-convex functions \cite{zhou2017stochastic,hinder2023near}.
For $\DRSRho = \rho \I$ it reduces to weak MVI, which was first introduced in \cite{diakonikolas2021Efficient} in the context of the extragradient method.
This class of nonmonotone problems emerges quite naturally in nonconvex-nonconcave minimax problems and nonconvex games, see for instance \cite[Ex. 6]{pethick2021Escaping}, \cite[Ex. 2.7]{evens2023convergence}.
It may also emerge in nonconvex optimization, as the following classical example demonstrates.

\begin{example}[Rosenbrock function]\label{ex:Rosenbrock}
    Let $b > 0$ and consider the function
    \begin{align*}
        f(x,y) \coloneqq x^2 + b(y-x^2)^2.
        \numberthis\label{eq:ex:Rosenbrock}
    \end{align*}
    The global minimum $(x^\star, y^\star) = (0, 0)$ of \eqref{eq:ex:Rosenbrock} is a $\Bigl(-\nicefrac{1}{16}\Bigr)$-weak Minty solution of $\nabla f$.
\end{example} 

To establish local convergence of \ref{eq:PPPA-intro}, we assume that \(T\) admits at least a single local \(\DRSRho\)-weak Minty solution. More specifically, we work under the following assumptions.

\begin{assumption}[local assumptions for \ref{eq:PPPA-intro}] \label{ass:PPPA-local}
    The operator \(T: \R^{n} \rightrightarrows \R^n\) and the linear mappings \(\M \in \sym{n}\) and \(\Lambda \in \sym{n}\) satisfy the following properties.
    \begin{enumeratass}
        \item \label{ass:PPPA-local:1} 
        There exists a nonempty, convex set \(\pazocal{S}^\star\subseteq \zer T\) and a symmetric, possibly indefinite matrix \(V \in \sym{n}\) such that \(T\) has $\DRSRho$\hyp{}oblique weak Minty solutions at \(\pazocal{S}^\star\) on a set $\pazocal U \subseteq \R^n \times \R^n$.
        \item \label{ass:PPPA-local:0} \(T\) is outer semicontinuous on $\pazocal U$.
        \item \label{ass:PPPA-local:2} \(\M\) and \(\Lambda\) are positive definite matrices satisfying \(\Lambda \M = \M \Lambda\) and
        \begin{equation}\label{eq:PPPA:eigcond-local}
            \Lambda
                \prec
            2(\I + \M^{\nicefrac12} \DRSRho \M^{\nicefrac12}).
        \end{equation}
        \item \label{ass:PPPA-local:0.5} There exists an $\varepsilon > 0$ such that for every 
        \[
            z \in \pazocal W_\varepsilon \coloneqq \set{z \in \R^n}[\dist_{\M \Lambda^{-1}}(z, \pazocal{S}^\star) \leq \varepsilon]
        \]
        there exists a \(\bar z\) satisfying
        \(
            \bigl(\bar z, P(z - \bar z)\bigr) \in \gph T \cap \pazocal U
        \).
        We drop the subscript and write \(\pazocal W\) whenever no ambiguity occurs.
    \end{enumeratass}
\end{assumption}

\begin{remark}
    \cref{ass:PPPA-local:2} imposes a restriction both on the choice of the preconditioner \(\M\) and the relaxation matrix \(\Lambda\).
    Since $\Lambda \succ 0$, a necessary condition for \eqref{eq:PPPA:eigcond-local} to hold is that \(\I + \M^{\nicefrac12} \DRSRho \M^{\nicefrac12} \succ 0\), which by positive definiteness of \(\M\) holds if and only if \(\M^{-1} + \DRSRho \succ 0\).
    This necessary condition is vacuously satisfied when $\DRSRho \succeq 0$, as is the case for monotone problems.
    Moreover, \eqref{eq:PPPA:eigcond-local} imposes an additional restriction on the relaxation matrix. For instance, when the relaxation matrix is a multiple of identity, i.e., $\Lambda = \ell\I$, \cref{ass:PPPA-local:2} reduces to the necessary condition
    \(
        \etamin
            {}\coloneqq{}
        \lambda_{\rm min}\bigl(\I + \M^{\nicefrac12} \DRSRho \M^{\nicefrac12}\bigr) > 0
    \)
    and the relaxation rule
    \(
        \ell \in (0, 2\etamin)
    \),
    showing that smaller relaxation parameters are required in the nonmonotone setting.
\end{remark}
\begin{remark}
    Note that the first step of \ref{eq:PPPA-intro} corresponds to finding a point $\bar{z}^k$ such that $\bigl(\bar{z}^k, P(z^k - \bar{z}^k)\bigr) \in \gph T$. In our local analysis, it is important to ensure that there exists such pairs within the local set \(\pazocal{U}\). 
    Provided that the current iterate $z^k$ is close enough to the solution set \(\pazocal{S}^\star\),
    this is ensured by \Cref{ass:PPPA-local:0.5}.
\end{remark}

\begin{figure}
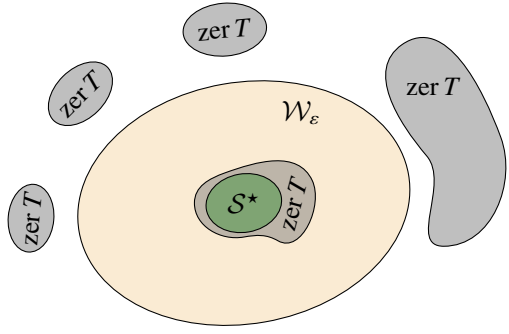

    \centering
    \includetikz{Local-assumptions/sets}
    \caption{
        Illustration of the subsets of $\R^n$ involved in \cref{ass:PPPA-local}.
        }
    \label{fig:local-assumptions}
\end{figure}
A visualization of the different subsets of $\R^n$ introduced in \cref{ass:PPPA-local} is provided in \Cref{fig:local-assumptions}.
Observe that in general our assumptions do not impose that the solution set \(\zer T\) is convex or even connected.
The situation changes when we focus on the class of maximally comonotone operators, which are a special case of operators satisfying the weak Minty condition.
For such operators, the set of zeros is necessarily convex as shown below. 

\begin{proposition}\label{prop:local-comon:convex-solutions}
    If an operator \(T : \R^n \rightrightarrows \R^n\) is maximally \(\DRSRho\)\hyp{}comonotone on a convex set
    \(
        \pazocal U \subseteq \R^n \times \R^n
    \), 
    then
    \(\zer T \cap \pazocal U\) is convex.
\end{proposition}

The following proposition demonstrates that our assumptions for \ref{eq:PPPA-intro} are satisfied for maximally comonotone operators, provided that the preconditioner \(\M \in \sym{n}\) and relaxation matrix \(\Lambda \in \sym{n}\) are selected accordingly.
The proof follows a similar reasoning as \cite[Lem. 1]{iusem2003Inexact} and is deferred to \Cref{proof:prop:local-comon}.


\begin{proposition}\label{prop:local-comon}
    Consider an operator \(T : \R^n \rightrightarrows \R^n\) for which $\zer T$ is nonempty,
    and let the preconditioner \(\M \in \sym{n}\) and relaxation matrix \(\Lambda \in \sym{n}\) satisfy \Cref{ass:PPPA-local:2}.
    Suppose that there exists a nonempty, convex set $\pazocal{S}^\star \subseteq \zer T$ and scalars $\delta_1 > 0$ and $\delta_2 > 0$ such that 
    \(T\) is maximally \(\DRSRho\)\hyp{}comonotone on
    \(
        \pazocal U_1 \times \pazocal U_2
    \), 
    where
    \[
        \pazocal U_1
            \coloneqq
        \set{z \in \R^n}[\dist_{\M \Lambda^{-1}}(z, \pazocal{S}^\star) \leq \delta_1]
        \quad\text{and}\quad
        \pazocal U_2
            \times
        \set{v \in \R^n}[\dist_{\M \Lambda^{-1}}(v, 0) \leq \delta_2].
    \]
    Then, \Cref{ass:PPPA-local} holds with
    \(
        \pazocal U
            {}={}
        \pazocal U_1 \times \pazocal U_2
    \)
    and
    \(
        \varepsilon
            =   
        \min
        \left(
            \frac{\delta_1}{\|\I - \Lambda^{-1}\| + \|\Lambda^{-1}\|},
            \frac{\delta_2}{
            \|\M\|\left(
                1 + \|\I - \Lambda^{-1}\| + \|\Lambda^{-1}\|
            \right)}
        \right).
    \)
    Moreover, if $\pazocal U_1 \times \pazocal U_2 = \R^n \times \R^n$, then this assumption holds globally.
\end{proposition}

Note that when \cref{ass:PPPA-local} holds globally, i.e., for 
\(\pazocal W = \R^n\) and \(\pazocal U = \R^n \times \R^n\),
then \Cref{ass:PPPA-local:1} and \Cref{ass:PPPA-local:0} simply reduce to assuming that $T$ has global $\DRSRho$-weak Minty solutions and is outer semicontinuous,
while \Cref{ass:PPPA-local:0.5} then corresponds to the preconditioned resolvent $(P + T)^{-1} \circ P$ having full domain. We discuss the global setting at the end of this section (see \cref{ass:PPPA-global}).

In our upcoming convergence analysis, we interpret each iteration of the preconditioned proximal point algorithm as a projection onto a certain halfspace. 
This perspective was previously employed in \cite{solodov1996modified,solodov1999hybrid,konnov1997Class} in the maximally monotone setting.
Specifically, for $\M = \gamma^{-1}\I$ and $\Lambda = \ell\I$, \cite{solodov1999hybrid} considers the halfspace
\begin{equation}
    \mathcal{D}_{z^k, \bar z^k} \coloneqq \set{r \in \R^n}[
        \langle \gamma^{-1} (z^k - \bar z^k), \bar z^k - r \rangle
        \ge 0].
\end{equation}
This halfspace was later generalized in \cite{evens2023convergence} to the form of \eqref{eq:hyperplane:Lambda} to account for symmetric positive (semi)definite preconditioners $\M$ and for operators $T$ that are not necessarily monotone but admit global weak Minty solutions.
In this work, we build on this idea by extending it to the weaker local setting and incorporating matrix relaxation $\Lambda$.
In particular, we show in the following lemma that this particular halfspace continues to separate the current iterate from a subset $\pazocal S^\star$ of solutions that satisfy the local oblique weak Minty assumption, provided that the current iterate $z$ is close enough to $\pazocal S^\star$.


\begin{lemma}[halfspace interpretation]\label{lem:hyperplane:Lambda}
    Suppose that \cref{ass:PPPA-local} holds with \(\pazocal{U}\) and \(\pazocal W\) denoting the sets therein.
    Let \(z \in \pazocal W\) and consider a point \(\bar z\) such that 
    \(
        \bigl(\bar z, P(z - \bar z)\bigr) \in \gph T \cap \pazocal U,
    \)
    which exists by \cref{ass:PPPA-local:0.5}.
    Let \(z^+ \coloneqq z + \Lambda (\bar z - z)\) and define the halfspace
    \begin{equation}\label{eq:hyperplane:Lambda}
        \mathcal{D}_{z, \bar z} \coloneqq \set{r \in \R^n}[
            \tinprod{\M(z - \bar z)}{\bar z - r}
            \ge \qindef{\M(z - \bar z)}{\DRSRho}].
    \end{equation}
    Then, the following hold.
    \begin{enumerate}
        \item\label{half:QS} The set \(\pazocal{S}^\star\) is a subset of \(\mathcal{D}_{z, \bar z}\).
        \item\label{half:xbar} If \(z \in \mathcal{D}_{z, \bar z}\) then
        \(z = \bar z\), and in particular \(\bar z \in \pazocal W \cap \zer T\).
        \item\label{half:xbar:2} 
        If \(z \notin \zer T\), then \(\|{z - \bar z}\|^2_{\M \Lambda} \neq 0\) and 
        \(
            \alpha \geq \bar\alpha > \nicefrac12
        \),
        where
        \begin{align*}
            \alpha 
                {}\coloneqq{}
            \frac{
                \|z-\bar{z}\|_{\M + \M\DRSRho\M}^2
            }
            {
                \|{z - \bar z}\|^2_{\M \Lambda}
            }
            \quad
            \text{and}
            \quad
            \bar \alpha
                {}\coloneqq{}
            \lambda_{\rm min}\left(\left(\I + \M^{\nicefrac12}\DRSRho\M^{\nicefrac12}\right)\Lambda^{-1}\right).
            \numberthis\label{eq:half:alpha}
        \end{align*}
        \item\label{half:update} 
        $
            z^+
                {}={}
            (1 - \nicefrac{1}{\alpha})z + \nicefrac{1}{\alpha}\proj^{\M \Lambda}_{\mathcal{D}_{z, \bar z}}(z)
        $, 
        and in particular for all \(z^\star \in \pazocal{S}^\star\) 
        \begin{align*}
            \|z^+-z^\star\|_{\M \Lambda^{-1}}^2
                {}\leq{}
            \|z-z^\star\|_{\M \Lambda^{-1}}^2 -
            (2\bar \alpha-1)
            \|z-\bar{z}\|_{\M \Lambda}^2.
            \numberthis\label{eq:half:update:descent}
        \end{align*}
    \end{enumerate}

    \begin{proof}
        \begin{proofitemize}
            \item \ref{half:QS}: 
            Since
            \(
                \bigl(\bar z, P(z - \bar z)\bigr) \in \gph T \cap \pazocal U,
            \)
            it holds that
            \begin{equation}
                \label{eq:PPPA-proof-vbar-incl}
                \M(z - \bar z) \in T \bar z.
            \end{equation}
            Thus, owing to the $\DRSRho$-oblique weak Minty assumption on \(\pazocal U\) (\cref{ass:PPPA-local:1}), it holds for all $z^\star \in \pazocal{S}^\star$ that
            \begin{equation}\label{eq:PPPA-VWMVI-application}
                \tinprod{\M(z - \bar z)}{\bar z - z^\star} \ge 
                \qindef{\M(z - \bar z)}{\DRSRho},
            \end{equation}
            showing that \(\pazocal{S}^\star \subseteq \mathcal{D}_{z, \bar z}\). 
            \item \ref{half:xbar}:
            Suppose that \(z\in\mathcal{D}_{z, \bar z}\)
            and
            define the shorthand notation \(\tilde{z} \coloneqq z - \bar z\).
            By definition of \(\mathcal{D}_{z, \bar z}\), it holds that
            \def\sqrtlam{X^{\nicefrac12}}
            \begin{align*}
                0 \geq \nrm{\tilde{z}}_{\M}^2 + \qindef{\M\tilde{z}}{\DRSRho}
                    {}={}&
                \tinprod{\M\tilde{z}}{(\I + \DRSRho\M)\tilde{z}}
                    {}={}
                \|\tilde{z}\|^2_{P+PVP}.
                \numberthis \label{eq:PPPA-ineq-posdef-X}
            \end{align*}
            Since
            \(\I + \M^{\nicefrac{1}{2}}\DRSRho \M^{\nicefrac{1}{2}}\) is positive definite owing to \cref{ass:PPPA-local:2}, this inequality implies that \(\tilde{z} = 0\), so that \(\bar z \in \zer T\) by \eqref{eq:PPPA-proof-vbar-incl}.
            
            \item \ref{half:xbar:2}:
            If \(\bar{z} \notin \zer T\), then inclusion \eqref{eq:PPPA-proof-vbar-incl} implies that \(\M(z - \bar z) \neq 0\).
            Since \(\M \succ 0\) and \( \Lambda \succ 0\) by assumption, it holds that \(\M\Lambda = \M^{\nicefrac{-1}{2}} \Lambda \M^{\nicefrac{-1}{2}} \succ 0\) and thus
            \(
                \|{z - \bar z}\|^2_{\M  \Lambda} 
                    {}>{}
                0
            \),
            ensuring that $\alpha$ as defined in \eqref{eq:half:alpha} is finite-valued.
            In particular, since $\alpha$ can be written as a Rayleigh quotient \cite[Thm. 4.2.2]{horn2012matrix} of \(\Lambda^{-1} + \Lambda^{\nicefrac{-1}2}\M^{\nicefrac12}\DRSRho\M^{\nicefrac12}\Lambda^{\nicefrac{-1}2}\) and \(\Lambda^{\nicefrac12} \tilde{z}\), it follows that
            \[
                \alpha
                    {}={}
                \frac{
                    \langle
                        \Lambda^{\nicefrac12}\M^{\nicefrac12}(z - \bar z),(\Lambda^{-1} + \Lambda^{\nicefrac{-1}2}\M^{\nicefrac12}\DRSRho\M^{\nicefrac12}\Lambda^{\nicefrac{-1}2}) \Lambda^{\nicefrac12} \M^{\nicefrac12}(z - \bar z)
                    \rangle
                }
                {
                    \|\Lambda^{\nicefrac12}\M^{\nicefrac12}(z - \bar z)\|^2
                }
                    {}\geq{}
                \lambda_{\rm min}(\Lambda^{-1} + \Lambda^{\nicefrac{-1}2}\M^{\nicefrac12}\DRSRho\M^{\nicefrac12}\Lambda^{\nicefrac{-1}2})
                    {}={}
                \bar \alpha,
            \]
            where the final equality holds by similarity transformation using the similarity matrix $\Lambda^{\nicefrac12}$ \cite[Cor. 1.3.4]{horn2012matrix}.
            The claim that $\bar \alpha > \nicefrac12$ follows directly by \cref{ass:PPPA-local:2} since the condition
            $2(\I + \M^{\nicefrac12}\DRSRho\M^{\nicefrac12}) \succ \Lambda$
            is equivalent to
            $\Lambda^{-1} + \Lambda^{\nicefrac{-1}2}\M^{\nicefrac12}\DRSRho\M^{\nicefrac12}\Lambda^{\nicefrac{-1}2} \succ \nicefrac12\I$.

            \item \ref{half:update}:
            Since \(P\) and \(\Lambda\) are positive definite, \( \mathcal{D}_{z, \bar z}\) may be equivalently represented as  
            \begin{equation*}
                \mathcal{D}_{z, \bar z} = \set{r \in \R^n}[
                    \tinprod{\Lambda(z - \bar z)}{\bar z - r}_{\M \Lambda^{-1}}
                    \geq
                    \qindef{\M(z - \bar z)}{\DRSRho}
                    ].
            \end{equation*}
            Consequently, for any $u \notin \mathcal{D}_{z, \bar z}$, it holds that \cite[Ex. 29.20]{bauschke2017Convex}
            \begin{align*}
                \proj^{\M \Lambda^{-1}}_{\mathcal{D}_{z, \bar z}}(u)
                    {}={}&
                u 
                    {}+{}  
                \frac{
                        \langle z-\bar{z},u - \bar{z}\rangle_{\M} + \qindef{\M(z-\bar{z})}{\DRSRho}
                    }
                    {
                        \|\Lambda(z - \bar z)\|_{\M \Lambda^{-1}}^2
                    }
                \bigl(\Lambda(\bar z - z)\bigr).
            \end{align*}
            Therefore,
            \(
                \proj^{\M \Lambda^{-1}}_{\mathcal{D}_{z, \bar z}}(z)
                    {}={}
                z
                    {}+{}
                \alpha \Lambda (\bar z - z)
            \)
            and it holds that
            \[
                z^+
                    {}={}
                z + \Lambda(\bar{z}-z)
                    {}={}
                (1 - \nicefrac{1}{\alpha})z + \nicefrac{1}{\alpha}\proj^{\M \Lambda^{-1}}_{\mathcal{D}_{z, \bar z}}(z)
                    {}\eqqcolon{}
                \tilde \proj(z).
            \]
            Owing to firm nonexpansiveness of $\proj^{\M \Lambda^{-1}}_{\mathcal{D}_{z, \bar z}}$ \cite[Prop. 4.16]{bauschke2017Convex}, the relaxed projection $\tilde \proj$ is $\nicefrac{1}{2\alpha}$-averaged
            in the space with inner product \(\langle \cdot,\cdot \rangle_{\M \Lambda^{-1}}\) \cite[Cor. 4.41]{bauschke2017Convex}.
            Note that for any $z^\star\in \pazocal{S}^\star$ it holds by \cref{half:QS} that $z^\star\in \mathcal{D}_{z, \bar z}$.
            Therefore, it follows from~\cite[Prop. 4.35(iii)]{bauschke2017Convex} for any $z^\star\in \pazocal{S}^\star$ that
            \begin{align}\label{eq:iFejer}
                \|z^+-z^\star\|_{\M \Lambda^{-1}}^2
                    {}={}
                \|\tilde \proj(z)-\tilde \proj(z^\star)\|_{\M \Lambda^{-1}}^2
                    &\leq
                \|z-z^\star\|_{\M \Lambda^{-1}}^2-\tfrac{1-\nicefrac{1}{2\alpha}}{\nicefrac{1}{2\alpha}}\|z^+-z\|_{\M \Lambda^{-1}}^2\nonumber\\
                &=\|z-z^\star\|_{\M \Lambda^{-1}}^2-(2\alpha-1)\|z-\bar{z}\|_{\M \Lambda}^2.
            \end{align}
            The claim is established by using that $\alpha \geq \bar \alpha$ owing to \Cref{half:xbar:2}.
            \qedhere
        \end{proofitemize}    
    \end{proof}
\end{lemma}

\cref{half:update} characterizes the update rule of \ref{eq:PPPA-intro} as an oblique projection onto a halfspace, which separates the current iterate from the subset of solutions \(\pazocal{S}^\star\) (as shown in \cref{half:QS}). 
Leveraging this result alongside the descent inequality in \eqref{eq:half:update:descent}, we can derive the following local convergence result for \ref{eq:PPPA-intro}, where the preconditioned resolvent evaluations of the operator \(T\) are restricted to the localization \(\pazocal{U}\) as in update rule \eqref{eq:PPPA-local}.
By construction, any sequence satisfying \eqref{eq:PPPA-local} also satisfies \eqref{eq:PPPA-intro}, and both update rules coincide when \(\pazocal{U} = \R^n \times \R^n\).

\begin{theorem}[local convergence of \ref{eq:PPPA-intro}]\label{thm:pppa-local}
    Suppose that \cref{ass:PPPA-local} holds with \(\pazocal{U}\) and \(\pazocal W\) denoting the sets therein.
    If $z^0 \in \pazocal W$, 
    then for any sequence $\seq{z^k, \bar z^k}$ generated by the update rule
    \begin{equation}\label{eq:PPPA-local}\tag{local-PPPA}
        \begin{cases}
            \text{find } \bar{z}^k \text{ such that } \bigl(\bar{z}^k, P(z^k - \bar{z}^k)\bigr) \in \gph T \cap \pazocal U\\
            z^{k+1}=z^k+\Lambda(\bar{z}^k-z^k)
        \end{cases}
    \end{equation}
    it holds that
    either a point $\bar z^k \in \pazocal W \cap \zer T$ is reached in a finite number of iterations or the following hold.
    \begin{enumerate}
        \item\label{it:pppa-local:v} 
        \(z^k\) belongs to \(\pazocal W\) for all \(k\in \N\)
        since
        \(
            \seq{\dist_{\M \Lambda^{-1}}(z^{k}, \pazocal{S}^\star)}
        \)
        is nonincreasing.
        \item\label{it:pppa-local:rate} 
        The sequences
        \(
            \seq{\dist_{\M \Lambda}^2\bigl(0, T(\bar z^k)\bigr)}
        \) 
        and
        \(
            \seq{\|z^k - \bar z^k\|_{\M \Lambda}}
        \)
        converge to zero with rate
        \begin{align*}
            \min_{k=0,1,\ldots,N} 
            \dist_{\M \Lambda}^2\bigl(0, T(\bar z^k)\bigr)
                {}\leq{}
            \min_{k=0,1,\ldots,N} 
            \|z^k - \bar z^k\|_{\M \Lambda}
                {}\leq{}
            \frac{\dist_{\M \Lambda^{-1}}^2(z^0, \pazocal{S}^\star)}{(N+1)(2\bar\alpha-1)},
            \numberthis\label{eq:it:pppa-local:rate}
        \end{align*}
        where \(\bar\alpha > \nicefrac12\) is defined as in equation \eqref{eq:half:alpha}.
        Additionally, if $T$ is $\DRSRho$\hyp{}comonotone on $\pazocal{U}$, then the sequence  
        \(
            \seq{\|z^k - \bar z^k\|_{\M \Lambda}}
        \) 
        is nonincreasing and it holds that
        \begin{align*}
            \dist_{\M \Lambda}^2\bigl(0, T(\bar z^N)\bigr)
                {}\leq{}
            \|z^N - \bar z^N\|_{\M \Lambda}
                {}\leq{}
            \frac{\dist_{\M \Lambda^{-1}}^2(z^0, \pazocal{S}^\star)}{(N+1)(2\bar\alpha-1)}.
            \numberthis\label{eq:it:pppa-local:rate:comon}
        \end{align*}
        \item\label{it:pppa-local:bounded} The sequences $(z^k)_{k\in\N}$, $(\bar{z}^k)_{k\in\N}$ are bounded and their limit points belong to $\pazocal W \cap \zer T$.
        \item\label{it:pppa-local:full} 
        If in \cref{ass:PPPA-local:1} the set
        $ 
            \pazocal{S}^\star
        $
        is equal to
        $ 
            \pazocal W \cap \zer T
        $,
        then the sequences
        $
            \seq{z^k}
        $,
        $
            \seq{\bar z^k}
        $
        converge
        to some element of $\pazocal W \cap \zer T$.
    \end{enumerate}
    \begin{proof}
        \begin{proofitemize}
            \item \ref{it:pppa-local:v}:
            Suppose for some \(k \geq 0\) that \(z^k \in \pazocal W\). 
            Then, \Cref{ass:PPPA-local:0.5} ensures that there exists a $\bar{z}^k$ satisfying update rule \eqref{eq:PPPA-local}.
            If $z^k \in \mathcal{D}_{z^k, \bar z^k}$ then by \Cref{half:xbar} the algorithm has reached a point $\bar{z}^k \in \zer T$.
            Otherwise, it follows from \Cref{half:update} for any $z^\star\in \pazocal{S}^\star$ that
            \begin{align}\label{eq:proof:pppa-local:iFejer}
                \|z^{k+1}-z^\star\|_{\M \Lambda^{-1}}^2
                    {}\leq{}
                \|z^k-z^\star\|_{\M \Lambda^{-1}}^2-(2\bar\alpha-1)\|z^k-\bar{z}^k\|_{\M \Lambda}^2.
            \end{align}
            Since $\bar\alpha > \nicefrac12$, this implies that
            \begin{align*}
                \dist_{\M \Lambda^{-1}}(z^{k+1}, \pazocal{S}^\star)
                    {}={}
                \inf_{z \in \pazocal{S}^\star} \left\{\|z^{k+1} - z\|_{\M \Lambda^{-1}}\right\}
                    {}\leq{}&
                \|z^{k+1}-\proj^{\M \Lambda^{-1}}_{\pazocal{S}^\star}(z^k, \pazocal{S}^\star)\|_{\M \Lambda^{-1}}\\
                    {}\overrel[\leq]{\eqref{eq:proof:pppa-local:iFejer}}{}&
                \|z^{k}-\proj^{\M \Lambda^{-1}}_{\pazocal{S}^\star}(z^{k}, \pazocal{S}^\star)\|_{\M \Lambda^{-1}}
                    {}={}
                \dist_{\M \Lambda^{-1}}(z^{k}, \pazocal{S}^\star),
            \end{align*}
            and therefore $z^{k+1} \in \pazocal W$. 
            Since we assumed that $z^0 \in \pazocal W$, it follows that either a point $\bar{z}^k \in \zer T$ is reached in a finite number of iterations, or that
            for all
            \(
                k \in \N
            \)
            it holds that
            \(
                z^k \in \pazocal W
            \),
            that
            \(
                \left(\bar z^k, \M(z^k - \bar z^k)\right)
                    \in 
                \gph T \cap \pazocal U
            \)
            and that \eqref{eq:proof:pppa-local:iFejer} holds for any $z^\star\in \pazocal{S}^\star$,
            establishing that $(z^k)_{k\in\N}$ is Fej\'er monotone with respect to $\pazocal{S}^\star$ \cite[Def. 5.1]{bauschke2017Convex}.

            \item \ref{it:pppa-local:rate}:
            By telescoping inequality \eqref{eq:proof:pppa-local:iFejer} for $k = 0, 1,\hdots,N$, we get that
            \(
                0
                    {}\leq{}
                \|z^{N+1}-z^\star\|_{\M \Lambda^{-1}}^2
                    {}\leq{}
                \|z^0-z^\star\|_{\M \Lambda^{-1}}^2
                    -
                \sum_{k=0}^{N}(2\bar\alpha-1)\|z^k-\bar{z}^k\|_{\M \Lambda}^2
            \) for any $z^\star\in \pazocal{S}^\star$.
            By rearranging and taking the square root, this implies that
            \begin{align*}
                \min_{k=0,1,\ldots,N}
                \|z^k-\bar{z}^k\|_{\M \Lambda}
                    {}\leq{}
                \frac{
                    \inf_{z^\star \in \pazocal{S}^\star} \left\{\|z^0-z^\star\|_{\M \Lambda^{-1}}\right\}
                }{
                    \sqrt{N+1}\sqrt{2\bar\alpha-1}
                }
                    {}={}
                \frac{
                    \dist_{\M \Lambda^{-1}}(z^0, \pazocal{S}^\star)
                }{
                    \sqrt{N+1}\sqrt{2\bar\alpha-1}
                }.
                \numberthis\label{eq:proof:it:pppa-local:rate}
            \end{align*}
            Then, the claimed rate from \eqref{eq:it:pppa-local:rate} follows by noting for all $k$ that
            \begin{align*}
                \dist_{\M^{-1} \Lambda}\bigl(0, T(\bar z^k)\bigr)
                    {}={}
                \min_{v \in T(\bar z^k)}
                \left\{
                    \|v\|_{\M^{-1} \Lambda}
                \right\}
                    {}\leq{}
                \|\M(z^k-\bar{z}^k)\|_{\M^{-1} \Lambda}
                    {}={}
                \|z^k-\bar{z}^k\|_{\M \Lambda},
                \numberthis\label{eq:proof:it:pppa-local:distance}
            \end{align*}
            where we used that
            \(
                \M(z^k - \bar z^k)
                    \in 
                T(\bar z^k)
            \).
            Define $\tilde z^k \coloneqq z^k - \bar z^k$ for all $k \in \N$ and note that
            \(
                \bigl(\bar{z}^k, \M\tilde z^k\bigr) \in \gph T \cap \pazocal U
            \)
            and
            \(
                \bigl(\bar{z}^{k+1}, \M\tilde z^{k+1}\bigr) \in \gph T \cap \pazocal U
            \).
            Consequently, if $T$ is $\DRSRho$\hyp{}comonotone on $\pazocal{U}$, then it holds that
            \(
                \inner{\bar z^k - \bar z^{k+1}, \M(\tilde z^k - \tilde z^{k+1})}
                    {}\geq{}
                \qindef{\M(\tilde z^k - \tilde z^{k+1})}{\DRSRho}.
            \)
            By reordering the terms and using the definition of $\tilde z^k$, this implies that
            \begin{align*}
                \inner{z^k - z^{k+1}, \M(\tilde z^k - \tilde z^{k+1})}
                    {}\geq{}&
                \qindef{\tilde z^k - \tilde z^{k+1}}{\M + \M \DRSRho \M}\\
                    {}={}&
                \qindef{\tilde z^k - \tilde z^{k+1}}{\M^{\nicefrac12}(\I + \M^{\nicefrac12} \DRSRho \M^{\nicefrac12})\M^{\nicefrac12}}
                    {}\overrel[\geq]{\eqref{eq:PPPA:eigcond-local}}{}
                \tfrac12\|\tilde z^k - \tilde z^{k+1}\|_{\M^{\nicefrac12}\Lambda\M^{\nicefrac12}}^2.
            \end{align*}
            Recall that \(\Lambda \M = \M \Lambda\) due to \Cref{ass:PPPA-local:2}.
            Plugging in \(z^{k+1}=z^k-\Lambda\tilde z^k\) from \eqref{eq:PPPA-local}, it follows that
            \(
                \inner{\tilde z^k, \M\Lambda(\tilde z^k - \tilde z^{k+1})}
                    {}\geq{}
                \tfrac12\|\tilde z^k - \tilde z^{k+1}\|_{\M\Lambda}^2.
            \)
            In turn, since
            $\|u\|_{\M\Lambda}^2 - \|v\|_{\M\Lambda}^2 = \inner{2\M\Lambda u, u-v} - \|u-v\|_{\M\Lambda}^2$
            for any $u, v \in \R^n$, this implies that
            \begin{align*}
                \|\tilde z^k\|_{\M\Lambda}^2 - \|\tilde z^{k+1}\|_{\M\Lambda}^2
                    {}={}
                \inner{2\M \Lambda \tilde z^k, \tilde z^k - \tilde z^{k+1}} - \|\tilde z^k - \tilde z^{k+1}\|_{\M\Lambda}^2
                    {}\geq{}
                0,
            \end{align*}
            showing that
            $\|\tilde{z}^k\|_{\M\Lambda}^2$
            is nonincreasing. Combining this with \eqref{eq:proof:it:pppa-local:rate} and \eqref{eq:proof:it:pppa-local:distance} yields the rate from \eqref{eq:it:pppa-local:rate:comon}.

            \item \ref{it:pppa-local:bounded}:
            By update rule~\eqref{eq:PPPA-local}, it follows that
            \(
                \bar v^k \coloneqq \M(z^k - \bar{z}^k) \in T \bar{z}^k    
            \)
            for all $k$.
            Since
            \(
                \seq{\|z^k - \bar z^k\|_{\M \Lambda}}
            \)
            converges to zero as shown in \Cref{it:pppa-local:rate} so does $\seq{\bar{v}^k}$, since
            \begin{align*}
                \|\bar{v}^k\|^2
                    {}={}
                \|\M(z^k-\bar{z}^k)\|^2
                    {}\leq{}
                \|\M^{\nicefrac12}\Lambda^{\nicefrac{-1}2}\|^2\, \|\Lambda^{\nicefrac12}\M^{\nicefrac12}(z^k-\bar{z}^k)\|^2
                    {}={}
                \nrm{\M \Lambda^{-1}}\, \|z^k-\bar{z}^k\|^2_{\M \Lambda}.
                \numberthis\label{eq:vw}
            \end{align*}
            Consequently, it follows from outer semicontinuity of $T$ on $\pazocal U$ (due to \Cref{ass:PPPA-local:0}) that any limit point of $(\bar{z}^k)_{k\in\N}$ belongs to $\zer T$.
            Moreover, it follows from \eqref{eq:proof:pppa-local:iFejer} that \(\seq{\|z^k - z^\star\|_{\M \Lambda^{-1}}}\) converges, and in particular that $(z^k)_{k\in\N}$ is bounded. In turn, using 
            that $\|z^k - \bar z^k\|_{\M \Lambda} \to 0$
            and the triangle inequality 
            \(
                \|\bar{z}^k\|_{\M \Lambda} \leq \|\bar{z}^k-z^k\|_{\M \Lambda}+\|z^k\|_{\M \Lambda},
            \)
            it follows that \(\seq{\bar{z}^k}\) is bounded and thus that it has at least one limit point.
            Take a subsequence $(\bar{z}^k)_{k\in K}$ converging to some limit point $z^\infty \in \zer T$.
            Since $(\bar{z}^k-z^k)_{k\in\N}$ converges to zero
            as shown in \Cref{it:pppa-local:v},
            we have that $(z^k)_{k\in K}$ also converges to the same limit point $z^\infty$.
            The claim is established by noting that $z^k \in \pazocal W$ for all $k$.

            \item \ref{it:pppa-local:full}:
            If
            $
                \pazocal{S}^\star
                    =
                \pazocal W \cap \zer T
            $,
            then the convergence of $(z^k)_{k\in\N}$ to some element of $\pazocal W \cap \zer T$ follows by~\cite[Thm. 5.5]{bauschke2017Convex}, and $(\bar z^k)_{k\in\N}$ converges to the same point owing to \Cref{it:pppa-local:v}.
            \qedhere
        \end{proofitemize}
    \end{proof}
\end{theorem}

Observe that our local analysis concerns the set $\pazocal W \cap \zer T$, i.e., the zeros of $T$ which are located within $\pazocal W$.
In that sense, the condition
$\pazocal{S}^\star = \pazocal W \cap \zer T$
as required in \cref{it:pppa-local:full}
simply states that all the zeros of $T$ within the set $\pazocal W$ are local weak Minty solutions (see also \cref{fig:local-assumptions}).

For completeness,  we conclude this section by providing the global counterpart of our results, considering the case where \(\pazocal W = \R^n\) and \(\pazocal U = \R^n \times \R^n\). In this setting, \cref{ass:PPPA-local} simplifies as follows.

\begin{assumption}[global assumptions for \ref{eq:PPPA-intro}] \label{ass:PPPA-global}
    The operator \(T: \R^{n} \rightrightarrows \R^n\) and the linear mappings \(\M \in \sym{n}\) and \(\Lambda \in \sym{n}\) satisfy the following properties.
    \begin{enumeratass}
        \item \label{ass:PPPA-global:1} 
        There exists a nonempty set \(\pazocal{S}^\star\subseteq \zer T\) and a symmetric, possibly indefinite matrix \(V \in \sym{n}\) such that \(T\) has $\DRSRho$\hyp{}oblique weak Minty solutions at \(\pazocal{S}^\star\).
        \item \label{ass:PPPA-global:0} \(T\) is outer semicontinuous.
        \item \label{ass:PPPA-global:2} \(\M\) and \(\Lambda\) are positive definite matrices satisfying \(\Lambda \M = \M \Lambda\) and
        \begin{equation}\label{eq:PPPA:eigcond}
            \Lambda
                \prec
            2(\I + \M^{\nicefrac12} \DRSRho \M^{\nicefrac12}).
        \end{equation}
        \item \label{ass:PPPA-global:0.5} The preconditioned resolvent $(P + T)^{-1} \circ P$ has full domain.
    \end{enumeratass}
\end{assumption}

Under these global conditions, the relaxed preconditioned proximal point algorithm is globally convergent, as summarized below. This result directly extends \cite[Thm. 2.4]{evens2023convergence} to incorporate matrix relaxation \(\Lambda\).

\begin{corollary}[global convergence of \ref{eq:PPPA-intro}]\label{thm:pppa}
    Suppose that \cref{ass:PPPA-global} holds. Consider a sequence $\seq{z^k, \bar z^k}$ generated by \ref{eq:PPPA-intro} starting from any $z^0 \in \R^n$.
    Then, either a point $\bar z^k \in \zer T$ is reached in a finite number of iterations or all the claims from \Cref{thm:pppa-local} hold with $\pazocal W = \R^n$.
\end{corollary}

\section{Local convergence of progressive decoupling+}\label{sec:progdec}

Recall that an equivalence relation between \ref{eq:progdec} and \ref{eq:PPPA-intro} was established in \cref{lem:equivalence:PPPA-progdec}.
In this section, we obtain local convergence results for \ref{eq:progdec} leveraging those for \ref{eq:PPPA-intro} from \cref{sec:local-pppa}, provided that the operator $S^X$ satisfies \cref{ass:PPPA-local}.
We will delve deeper into this idea and provide explicit assumptions on the operator $S$ for establishing convergence of \ref{eq:progdec}, rather than assumptions on the partial inverse $S^X$.
To achieve this, we first present the class of semimonotone operators, which is discussed in more detail in \cite{evens2023convergence, evens2023convergenceCP, quan2024scaled}.

\begin{definition}[semimonotonicity]\label{def:semimonotonicity}
    Let $\Mon,\Com\in\R^{n \times n}$ be symmetric (possibly indefinite) matrices. 
    An operator 
    $T : \R^n \rightrightarrows \R^n$ is said to be $(\Mon,\Com)$\hyp{}semimonotone at $(\other{x}, \other{y}) \in \gph T$ on a set $\pazocal U \ni (\other{x}, \other{y})$ if
    \begin{equation}\label{eq:obliquequasisemimonotonicity:matrix}
        \inner*{x - \other{x}, y - \other{y}} \geq \qindef{x - \other{x}}{\Mon} + \qindef{y - \other{y}}{\Com}, \qquad \text{for all } (x, y) \in \gph T \cap \pazocal U,
    \end{equation}
    where ${\displaystyle \qindef{\cdot}{X} \dfn \inner{\cdot,\cdot}_X}$ for any symmetric matrix $X\in\R^{n \times n}$.

    An operator \(T\) is said to be $(\Mon,\Com)$\hyp{}semimonotone on $\pazocal U$ if it is $(\Mon,\Com)$\hyp{}semimonotone at all $(\other{x}, \other{y})\in \gph T \cap \pazocal U$.
    It is said to be maximally \((\Mon,\Com)\)-semimonotone on $\pazocal U$ if its graph is not strictly contained in the graph of another \((\Mon,\Com)\)-semimonotone operator on $\pazocal U$. 

	Throughout, whenever $\pazocal U = \R^n \times \R^n$, the set $\pazocal U$ is omitted, 
	and whenever
    \(\Mon = \mon \I_n\) and \(\Com = \com \I_n\) where $\mon, \com \in \R$, 
    the prefix $(\Mon, \Com)$ is replaced by $(\mon,\com)$ and condition \eqref{eq:obliquequasisemimonotonicity:matrix} reduces to
    \begin{equation}\label{eq:quasisemimonotonicity}
        \inner*{x - \other{x}, y - \other{y}} \geq \mon \nrm{x - \other{x}}^2 + \com \nrm{y - \other{y}}^2, \qquad \text{for all } (x, y) \in \gph T \cap \pazocal U.
    \end{equation}
\end{definition}

Several special cases of semimonotonicity are particularly relevant in the context of this paper and are listed below.  
For additional connections to other classes of operators, we refer to \cite[Rem. 4.2]{evens2023convergence}.
\begin{remark}[relationship with other types of operators]\label{rem:relationship:WMVI}
    \begin{enumerate}
		\item \((\mon, 0)\)-semimonotonicity is equivalent to \(\mon\)-monotonicity, which is also known as \(|\mon|\)-\emph{hypo}-monotonicity when \(\mon < 0\), as monotonicity when \(\mon = 0\), and as \emph{strong} monotonicity when \(\mon > 0\).
		\((0, \com)\)-semimonotonicity on the other hand is equivalent to \(\com\)-comonotonicity, which is also referred to as \(|\com|\)-co\emph{hypo}-monotonicity when \(\com < 0\), and as \(\com\)-cocoercivity when \(\com > 0\).
		\item 
		An operator has \(\DRSRho\)-oblique weak Minty solutions at \(\pazocal{S}^\star \subseteq \zer T\) on \(\pazocal U\), as required in \cref{ass:PPPA-local:1},  
		if and only if it is \((0, \DRSRho)\)-semimonotone at every \((x^\star, 0) \in \graph T\) on \(\pazocal U\) for all \(x^\star \in \pazocal{S}^\star\).  
		Similarly, maximal \(\DRSRho\)-comonotonicity on \(\pazocal U\) and maximal \((0, \DRSRho)\)-semimonotonicity on \(\pazocal U\) describe identical conditions (see \cref{prop:local-comon}).
		\item Finally, an operator \(S\) is said to be \emph{elicitable monotone} at a level \(e \geq 0\) if \(S + e \proj_{X^{\bot}}\) is maximally monotone \cite[Def. 1]{rockafellar2019Progressive}, which is equivalent to maximal \((-e \proj_{X^{\bot}}, 0)\)\hyp{}semimonotonicity of \(S\).
        \qedhere
    \end{enumerate}
\end{remark}

In \cite{rockafellar2019Progressive}, the convergence of standard progressive decoupling was established in the elicitable monotone setting. However, many problems in optimization and variational analysis fall outside this particular framework (see e.g. \cref{it:ex:linear-system:2}).
Therefore, in our work we will focus on the case where $S$ belongs to the more general class of (locally) $(\mon \proj_{X^\bot}, \com \proj_X)$-semimonotone operators.
We proceed with several illustrative examples.
The first example is an equivalent formulation of the Rosenbrock function, which is shown to be semimonotone at its global minimum (see also \cref{ex:Rosenbrock}).
The second example involves finding the local minimum of an (inverted) double-well potential, which is shown to be only locally semimonotone around this local minimum.
The third example considers a matrix splitting problem, which is shown to be semimonotone due to its linearity.
The proofs are deferred to \Cref{proof:ex:Rosenbrock:progdec}.

\begin{example}[Rosenbrock function]\label{ex:Rosenbrock:progdec}
    Let $b > 0$ and consider the minimization problem
    \begin{align*}
        \minimize_{x \in \R^3}
        \quad
        f(x) \coloneqq x_1x_2 + b(x_3-x_1^2)^2
        \quad
        \stt
        \quad
        x \in X \coloneqq \set{x \in \R^3}[x_1 = x_2].
        \numberthis\label{eq:ex:Rosenbrock:progdec}
    \end{align*}
    Then, $\nabla f$ is 
    $\Bigl(-\nicefrac94 \proj_{X^\bot}, -\nicefrac{1}{4}\proj_X\Bigr)$-semimonotone
    at $(x^\star, 0)$, where $x^\star = (0, 0, 0)$ is the global minimum of \eqref{eq:ex:Rosenbrock}.
\end{example}

\begin{example}[local minimum]\label{ex:double-well}
    Consider the minimization problem
    \begin{align*}
        \minimize_{x \in \R^2}
        \quad
        f(x) \coloneqq \tfrac12 x_1^2 + \tfrac12x_2^2 - \tfrac14x_1^2x_2^2
        \quad
        \stt
        \quad
        x \in X \coloneqq \set{x \in \R^2}[x_1 = x_2].
        \numberthis\label{eq:example:double-well}
    \end{align*}
    This optimization problem is unbounded, but has a local minimum at $x^\star = (0, 0)$, and the operator $\nabla f$ is $(0, \proj_X)$-semimonotone at $(x^\star, 0)$ on $\pazocal U \coloneqq \set{x \in \R^2}[\|x - x^\star\| \leq 2] \times \R^2$.
\end{example}

\begin{example}[matrix splitting]\label{ex:consensus:cond}
    Consider finding an $x \in \R^n$ such that
    \(
        \sum_{i=1}^N A_i x = b
    \),
    where \(A_i \in \R^{n \times n}\) for $i = 1, \hdots, N$ and $b \in \R^n$.
    Defining the matrix
    \(A \coloneqq \blkdiag(A_1, \hdots, A_N)\), the operator
    \(
        S(x) \coloneqq A(x)
        - 
        \begin{psmallmatrix}
            0_{(N-1)n}\\
            b
        \end{psmallmatrix}
    \)
    and the consensus set 
    \(
        X \coloneqq \smash{\set{x \in \R^{Nn}}[x_1 = \hdots = x_N]}
    \),
    this is equivalent to solving linkage problem \eqref{prob:linkage-intro}.
    Let
    \(
        \mon \geq 0
    \)
    and
    \(
        \com \leq 0
    \).
    If the matrices $A_i$ are \((\mon, \com)\)-semimonotone for all $i = 1, \hdots, N$, and if
    \begin{equation}\label{eq:ex:linear:StS0}
        A^\top11^\top A \succeq \nu A^\top A
    \end{equation}
    for some \(\nu >0\),
    then operator \(S\) is
    $(\mon \proj_{X^\bot}, \frac{N\com}{\nu} \proj_X)$\hyp{}semimonotone.
\end{example}

As a particular instance, note that if \(A_i^\top A_j = 0 \) for all \(i\neq j\), then \eqref{eq:ex:linear:StS0} holds with equality and \(\nu = 1\).

Some additional examples of $(\mon \proj_{X^\bot}, \com \proj_X)$-semimonotone mappings to which our upcoming theory applies will be presented in \cref{subsec:progdec:examples}, including a numerical matrix splitting problem.

Our primary interest in this class of operators stems from \Cref{cor:partialinverse:semi-WMVI}, which establishes an equivalence between the $(\mon \proj_{X^\bot}, \com \proj_X)$-semimonotonicity of $S$ and the existence of oblique weak Minty solutions for its partial inverse $S^X$.  
In that sense, this particular operator class naturally extends the weak MVI condition considered previously for \ref{eq:PPPA-intro}. This corollary is a direct consequence of \Cref{prop:partialinverse:S-SX}, which is deferred to the appendix.


\begin{corollary}[partial inverse] 
	\label{cor:partialinverse:semi-WMVI}
	Let \(\mon, \com \in \R\) and consider an operator \(S:\R^n\rightrightarrows\R^n\) and a closed linear subspace \(X\subseteq \R^n\).
	Suppose that $\link_X S$ is nonempty.
	Then, the following are equivalent.
	\begin{enumerate}
		\item\label{it:cor:partialinverse:semi} Operator $S$ is $\bigl(\mon\proj_{X^{\bot}}, \com\proj_{X}\bigr)$\hyp{}semimonotone at 
		$(x^\star, y^\star) \in \link_X S$
		on $\pazocal U$.
		\item\label{it:cor:partialinverse:WMVI} Operator \(S^X\) has a \(\bigl(\mon \proj_{X^\bot}+\com \proj_X \bigr)\)-oblique weak Minty solution at 
		$
			\pazocal{S}^\star 
				\coloneqq
			\set{x^\star + y^\star}
				\subseteq
			\zer S^X
		$
		on 
		\(
			L_X(\pazocal U).
		\)
	\end{enumerate}
\end{corollary}

Leveraging this equivalence result, we can directly translate the conditions on the partial inverse $S^X$ for \ref{eq:PPPA-intro} to equivalent conditions on the operator $S$ for \ref{eq:progdec}, leading to the following set of assumptions.

\begin{assumption}[local assumptions for \ref{eq:progdec}] \label{ass:progdec-local}
	The operator \(S: \R^{n} \rightrightarrows \R^n\) and the stepsizes \(\gamma, \lambda_x, \lambda_y \in \R\) satisfy the following properties.
	\begin{enumeratass}
		\item \label{ass:progdec-local:3} 
		There exist parameters \(\mon, \com \in \R\)
		satisfying
		\(
			[\mon]_{-}[\com]_{-} < 1
		\)
		such that
		$S$ is $(\mon \proj_{X^\bot}, \com \proj_X)$\hyp{}semimonotone at 
		a solution $(x^\star, y^\star) \in \link_X S$
		on a set $\pazocal U^{\rm link} \subseteq \R^n \times \R^n$.
		\item \label{ass:progdec-local:1} Operator \(S\) is outer semicontinuous on $\pazocal U^{\rm link}$. 
		\item \label{ass:progdec-local:stepsize} The stepsizes $\gamma$, $\lambda_x$ and $\lambda_y$ satisfy\footnote{Note that the stepsize ranges provided in \eqref{eq:thm:progdec-local:stepsize} are nonempty since \([\mon]_{-}[\com]_{-} < 1\).}
		\begin{align*}\numberthis\label{eq:thm:progdec-local:stepsize}
			\gamma \in \bigl([\mon]_-, \nicefrac{1}{[\com]_-}\bigr),
			\quad
			\lambda_x \in \bigl(0, 2(1 + \gamma\com)\bigr)
			\quad\text{and}\quad
			\lambda_y \in \bigl(0, 2(1 + \nicefrac\mon\gamma)\bigr).
		\end{align*}
		\item \label{ass:progdec-local:2} For the selected stepsizes, 
		there exists an $\varepsilon > 0$ such that for every
		\[
            (x,y) \in \pazocal{W}^{\rm link}_\varepsilon
				\coloneqq 
			\set{
				(x,y) \in X \times X^\bot
			}[
				\gamma\lambda_x^{-1}\|x - x^\star\|^2
					+
				\gamma^{-1}\lambda_y^{-1}\|y - y^\star\|^2
					\leq
				\varepsilon^2
			],
        \]
        there exists a pair \((\bar x, \bar y)\) such that
        \begin{align*}
            \bigl(\bar x + \gamma^{-1}(y-\bar y), \bar y + \gamma(x - \bar x)\bigr) \in \gph S \cap \pazocal U^{\rm link}.
			\numberthis\label{eq:ass:progdec-local:2:cond}
		\end{align*}
		We drop the subscript and write \(\pazocal W^{\rm link}\) whenever no ambiguity occurs.
	\end{enumeratass}
\end{assumption}
The equivalence of our assumptions is formally shown in \cref{lem:equiv-ass}, whose proof involves showing that each individual item of \Cref{ass:progdec-local}
for \ref{eq:progdec} is indeed a direct translation of those from \Cref{ass:PPPA-local} for \ref{eq:PPPA-intro} when 
$T \coloneqq S^X$,
\(
	\M \coloneqq \gamma \proj_X + \gamma^{-1} \proj_{X^\bot}
\),
\(
	\Lambda \coloneqq \lambda_x \proj_X + \lambda_y \proj_{X^\bot}
\)
and
\(
	\DRSRho \coloneqq \mon \proj_{X^\bot} + \com\proj_X
\)
(recall the equivalence of \ref{eq:progdec} and \ref{eq:PPPA-intro} from \Cref{lem:equivalence:PPPA-progdec}).
\begin{lemma}\label{lem:equiv-ass}
	Consider an operator \(S:\R^n\rightrightarrows\R^n\) and a closed linear subspace \(X\subseteq \R^n\).
	Let \(\mon, \com, \gamma, \lambda_x, \lambda_y, \varepsilon \in \R\)
	and define
	$T \coloneqq S^X$,
	\(
		\M \coloneqq \gamma \proj_X + \gamma^{-1} \proj_{X^\bot}
	\),
	\(
		\Lambda \coloneqq \lambda_x \proj_X + \lambda_y \proj_{X^\bot}
	\)
	and
	\(
		\DRSRho \coloneqq \mon \proj_{X^\bot} + \com\proj_X
	\).
	Then, the following are equivalent.
	\begin{enumerate}
		\item
		\Cref{ass:progdec-local} 
		holds
		for $S, \mon, \com, \gamma, \lambda_x$ and $\lambda_y$
		with
		\(
			(x^\star, y^\star) \in \link_X S
		\)
		and
		\(
			\pazocal U^{\rm link} \subseteq \R^n \times \R^n
		\).
		\item 
		\Cref{ass:PPPA-local}
		holds 
		for $T, \M, \Lambda$ and $\DRSRho$
		with
		\(
			\pazocal{S}^\star
				{}={}
			\set{x^\star + y^\star}
				{}\subseteq{}
			\zer T
		\)
		and
		\ifsiam
\else
			\(
				\pazocal{U}
					{}={}
				L_X(\pazocal U^{\rm link})
				\subseteq \R^n \times \R^n
			\).
		\fi
	\end{enumerate}
	\begin{proof}
		Owing to
		\Cref{cor:partialinverse:semi-WMVI}, it follows that \cref{ass:progdec-local:3} holds 
		if and only if
		\(T = S^X\) has 
		\(
			\DRSRho
				{}={} 
			\bigl(
				\mon \proj_{X^\bot} + \com\proj_X
			\bigr)
		\)\hyp{}oblique weak Minty solutions at
		$
			\pazocal{S}^\star = \set{x^\star + y^\star} \subseteq \zer S^X
		$
		on 
		\(
			\pazocal{U}
				{}={}
			L_X(\pazocal U^{\rm link})
		\),
		i.e., if and only if \cref{ass:PPPA-local:1} holds.
		By definition of outer semicontinuity and the partial inverse, $S$ is outer semicontinuous on $\pazocal U^{\rm link}$ if and only if $S^X$ is outer semicontinuous on 
		$
			\pazocal{U}
				{}={}
			L_X(\pazocal U^{\rm link})
		$,
		establishing the equivalence between \Cref{ass:PPPA-local:0} and \cref{ass:progdec-local:1}.
		Applying the partial inverse to \eqref{eq:ass:progdec-local:2:cond} and using that \(L_X(\gph S) = \gph S^X\),
		it follows that
		\cref{ass:progdec-local:2} holds if and only if
		there exists an $\varepsilon > 0$ such that for every 
		\begin{align*}
            (x,y)
				{}\in{}&
			\set{
				(x,y) \in X \times X^\bot
			}[
				\gamma\lambda_x^{-1}\|x - x^\star\|^2
					+
				\gamma^{-1}\lambda_y^{-1}\|y - y^\star\|^2
					\leq
				\varepsilon^2
			]\\
				{}={}&
			\set{
				(x,y) \in X \times X^\bot
			}[
				\dist_{
						\left(
							\gamma\lambda_x^{-1}\proj_X + \gamma^{-1}\lambda_y^{-1}\proj_{X^\bot}
						\right)
					}
					\left(
						x + y, \pazocal{S}^\star
					\right) \leq \varepsilon
			],
        \end{align*}
		there exists a pair \((\bar x, \bar y)\) such that
		\(
			\bigl(\bar x + \bar y, \gamma(x - \bar x) + \gamma^{-1}(y-\bar y)\bigr) \in \gph S^X \cap L_X(\pazocal U^{\rm link}).
		\)
		Defining \(z \coloneqq x + y\) and \(\bar z \coloneqq \bar x + \bar y\) and noting that $x \in X, \bar x \in X, y \in X^{\bot}$ and $\bar y \in X^{\bot}$, this is equivalent to stating that for every
		\[
			z \in \set{z \in \R^n}[\dist_{\M \Lambda^{-1}}(z, \pazocal{S}^\star) \leq \varepsilon],
		\]
		there exists a \(\bar z\) such that
		\(
			\bigl(\bar z, \M(z - \bar z)\bigr) \in \gph S^X \cap L_X(\pazocal U^{\rm link})
		\),
		matching \Cref{ass:PPPA-local:0.5}.

		It only remains to show that \cref{ass:PPPA-local:2} holds if and only if \eqref{eq:thm:progdec-local:stepsize} holds.
		First, observe that the preconditioner $\M$ and relaxation matrix $\Lambda$ are positive definite if and only if $\gamma > 0$, $\lambda_x > 0$ and $\lambda_y > 0$. 
		Moreover, using that $\I = \proj_X + \proj_{X^\bot}$ and plugging in $\M, \Lambda$ and $\DRSRho$, condition \eqref{eq:PPPA:eigcond} from \cref{ass:PPPA-local:2} is equivalent to having
		\[
			\lambda_x \proj_X + \lambda_y \proj_{X^\bot}
				\prec
			2(1 + \gamma \com) \proj_X + 2(1 + \nicefrac{\mon}{\gamma}) \proj_{X^\bot}.
		\]
		Therefore, indeed \cref{ass:PPPA-local:2} holds if and only if 
		\eqref{eq:thm:progdec-local:stepsize}
		is satisfied for $\gamma, \lambda_x$ and $\lambda_y$.
	\end{proof}
\end{lemma}

Due to its equivalence with \Cref{ass:PPPA-local} for \ref{eq:PPPA-intro}, \cref{ass:progdec-local} provides a set of mild conditions under which \ref{eq:progdec} converges locally.
For instance, the solution set $\link_X S$ of linkage problem \eqref{prob:linkage-intro} need not be convex or even connected, as visualized previously for \ref{eq:PPPA-intro} in \Cref{fig:local-assumptions}.

A notable special case of \Cref{ass:progdec-local} arises when \( S \) is maximally  
\(\bigl(\mon\proj_{X^{\bot}}, \com\proj_{X}\bigr)\)\hyp{}semimonotone on \( \pazocal{U} \),  
or equivalently, when the partial inverse of \( S \) is maximally  
\(\bigl(\mon \proj_{X^\bot}+\com \proj_X \bigr)\)\hyp{}comonotone on  
\( L_X(\pazocal{U}) \) (see \Cref{prop:partialinverse:S-SX}).  
The following proposition provides sufficient conditions for a maximally semimonotone operator \( S \)  
on \( \pazocal{U}^{\rm link}_1 \times \pazocal{U}^{\rm link}_2 \) to satisfy \Cref{ass:progdec-local},  
where \( \pazocal{U}^{\rm link}_1 \) and \( \pazocal{U}^{\rm link}_2 \) are balls as defined in \eqref{eq:prop:local-semi:balls}.  
To avoid explicitly computing \( L_X(\pazocal{U}^{\rm link}_1 \times \pazocal{U}^{\rm link}_2) \) in the proof,  
the key idea is to first construct subsets of \( X \) and \( X^\bot \) whose Minkowski sum fits within  
\( \pazocal{U}^{\rm link}_1 \) and \( \pazocal{U}^{\rm link}_2 \), and then apply the partial inverse \( L_X \)  
to these structured Minkowski sums instead. The proof is deferred to \Cref{proof:prop:local-semi}.


\begin{proposition}\label{prop:local-semi}
	Consider an operator \(S:\R^n\rightrightarrows\R^n\) and a closed linear subspace \(X\subseteq \R^n\).
	Let \(\mon, \com \in \R\) satisfy
	\(
		[\mon]_{-}[\com]_{-} < 1
	\)
	let $\gamma$, $\lambda_x$ and $\lambda_y$ satisfy
	\eqref{eq:thm:progdec-local:stepsize}
	and suppose that $\link_X S$ is nonempty.
	Suppose that there exists a point $(x^\star, y^\star) \in \link_X S$ and a scalar $\delta > 0$ such that 
    \(S\) is maximally $(\mon \proj_{X^\bot}, \com \proj_X)$\hyp{}semimonotone on
    \(
        \pazocal{U}^{\rm link}_1 \times \pazocal{U}^{\rm link}_2
    \), 
    where
    \begin{align*}
        \pazocal{U}^{\rm link}_1
            \coloneqq
        \set{z \in \R^n}[\|z - x^\star\| \leq \delta]
        \quad\text{and}\quad
        \pazocal{U}^{\rm link}_2
            \coloneqq
        \set{z \in \R^n}[\|z - y^\star\| \leq \delta].
		\numberthis\label{eq:prop:local-semi:balls}
	\end{align*}
    Then, \Cref{ass:progdec-local} holds with
	\(
		\pazocal{U}^{\rm link}
			=
		\pazocal{U}^{\rm link}_1 \times \pazocal{U}^{\rm link}_2	
	\)
	and
	\(
		\varepsilon
			=   
		\frac{
			\delta\sqrt{\min\bigl\{\gamma\lambda_x^{-1}, (\gamma\lambda_y)^{-1}\bigr\}}			
			\min(\gamma, \gamma^{-1})
		}{
			\sqrt{2}\left(
				1 - \min\{\lambda_x^{-1}, \lambda_y^{-1}\} + \max\{\lambda_x^{-1}, \lambda_y^{-1}\}
			\right)
		}.
	\)
    Moreover, if $\pazocal{U}^{\rm link}_1 \times \pazocal{U}^{\rm link}_2 = \R^n \times \R^n$, then this assumption holds globally.
\end{proposition}

Having established that \ref{eq:progdec} can be cast as an instance of \ref{eq:PPPA-intro} (see \Cref{lem:equivalence:PPPA-progdec}), and that our local assumptions for \ref{eq:progdec} directly translate from those in \cref{ass:PPPA-local} for \ref{eq:PPPA-intro} (see \cref{lem:equiv-ass}), we obtain the following local convergence result for \ref{eq:progdec} as an immediate consequence of \cref{thm:pppa-local} for \ref{eq:PPPA-intro}, where we restrict the resolvent evaluations of the operator \(S\) to the localization \(\pazocal{U}^{\rm link}\).

\begin{theorem}[local convergence of \ref{eq:progdec}]\label{thm:progdec-local}
	Suppose that \Cref{ass:progdec-local} holds with \(\pazocal{U}^{\rm link}\) and \(\pazocal{W}^{\rm link}\) denoting the sets therein.
	If $(x^0, y^0) \in \pazocal{W}^{\rm link}$, then for any 
	\ifsiam
	\else
		sequence
	\fi
	$(\bar{x}^k, \bar{y}^k, x^k, y^k)_{k\in\N}$ generated by the update rule
	\begin{equation*}
		\begin{cases}
			\text{find } \bar{x}^k \text{ and } \bar{y}^k \text{ such that }
			\bigl(\bar{x}^k + \gamma^{-1}(y^k-\bar{y}^k), \bar{y}^k + \gamma(x^k - \bar{x}^k)\bigr) \in \gph S \cap \pazocal U^{\rm link}\\
			x^{k+1}
				{}={}
			x^k + \lambda_x(\bar{x}^k - x^k)\\
			y^{k+1}
				{}={}
			y^k + \lambda_y(\bar{y}^k - y^k)
		\end{cases}
	\end{equation*}
	it holds that either a point $(\bar{x}^k, \bar{y}^k) \in \pazocal{W}^{\rm link} \cap \link_X S$ is reached in a finite number of iterations or the following hold.
	\begin{enumerate}
		\item\label{it:progdec-local:v}
		The pair $(x^k, y^k)$ belongs to $\pazocal{W}^{\rm link}$
        for all \(k\in \N\) since
		\(	
			\seq{		
				\gamma\lambda_x^{-1}\|x - x^\star\|^2
					+
				\gamma^{-1}\lambda_y^{-1}\|y - y^\star\|^2
			}
		\)
		is nonincreasing.
		\item\label{it:progdec-local:rate}
		Let 
		\(
			\bar \alpha \coloneqq \min\Bigl(\tfrac{1 + \gamma \com}{\lambda_x}, \tfrac{1 + \nicefrac{\mon}{\gamma}}{\lambda_y}\Bigr) > \nicefrac12
		\).
		The sequence
		\(
			\seq{\gamma\lambda_x\nrm{\bar{x}^k - x^k}^2 + \gamma^{-1}\lambda_y\nrm{\bar{y}^k - y^k}^2}
		\)
		satisfies
        \begin{align*}
            \min_{k=0,1,\ldots,N} 
            \gamma\lambda_x\nrm{\bar{x}^k - x^k}^2 + \gamma^{-1}\lambda_y\nrm{\bar{y}^k - y^k}^2
                {}\leq{}
            \frac{
				\gamma\lambda_x^{-1}\nrm{x^0-x^\star}^2 + \gamma^{-1}\lambda_y^{-1}\nrm{y^0-y^\star}^2
			}{
				(N+1)(2\bar\alpha-1)
			}.
			\numberthis\label{eq:it:progdec-local:rate}
        \end{align*}
		Additionally, if \(S\) is $(\mon \proj_{X^\bot}, \com \proj_X)$\hyp{}semimonotone on $\pazocal U^{\rm link}$, then
		this sequence is nonincreasing and the minimum in \eqref{eq:it:progdec-local:rate} is therefore attained for $k = N$.
		\item\label{it:progdec-local:bounded} The sequences $\seq{x^k, y^k}$ and $\seq{\bar x^k, \bar y^k}$ are bounded and their limit points belong to $\pazocal{W}^{\rm link} \cap \link_X S$.
		\item\label{it:progdec-local:full} 
		If in \cref{ass:progdec-local:3} the set
		$ 
			\pazocal{S}^\star
		$
		is equal to
		$ 
			\pazocal{W}^{\rm link} \cap \link_X S
		$,
		then the sequences
		$\seq{x^k, y^k}$,
		$\seq{\bar x^k, \bar y^k}$
		converge to some element of $\pazocal{W}^{\rm link} \cap \link_X S$.
	\end{enumerate}
	\begin{proof}
		Owing to \Cref{lem:equivalence:PPPA-progdec}, update rule \ref{eq:progdec} can be cast as an instance of \ref{eq:PPPA-intro}.
		Moreover, \cref{lem:equiv-ass} ensures that the operator $T = S^X$, 
		preconditioner
		$
			\M = \gamma \proj_X + \gamma^{-1} \proj_{X^\bot}
		$
		and relaxation matrix
		$
			\Lambda = \lambda_x \proj_X + \lambda_y \proj_{X^\bot}
		$
		satisfy \cref{ass:PPPA-local}.
		Consequently, all the assertions follow from those of \Cref{thm:pppa-local} applied $0 \in S^X(z)$. 
		In particular, the claims holds using the relations
		\(L_X(z^k, 0) = (\proj_X(z^k), \proj_{X^\bot}(z^k)) = (x^k, y^k)\),
		\(L_X(\bar z^k, 0) = (\proj_X(\bar z^k), \proj_{X^\bot}(\bar z^k)) = (\bar x^k, \bar y^k)\)
		and that
		\begin{align*}
			&{} z^\star \in{} 
			\set{z \in \R^n}[\dist_{\M \Lambda^{-1}}(z, \pazocal{S}^\star) \leq \varepsilon] \cap \zer S^X
			{}\;\Longleftrightarrow\;{}
			(\proj_X(z^\star), \proj_{X^\bot}(z^\star)) \in \pazocal{W}^{\rm link} \cap \link_X S.
			\qedhere
		\end{align*}
    \qedhere
	\end{proof}
\end{theorem}

When \cref{ass:progdec-local} holds globally, i.e., for \(\pazocal U = \R^n \times \R^n\) and \(\pazocal{W}^{\rm link} = X \times X^\bot\), our assumptions can be stated as follows.

\begin{assumption}[global assumptions for \ref{eq:progdec}] \label{ass:progdec-global}
	The operator \(S: \R^{n} \rightrightarrows \R^n\) and the stepsize \(\gamma \in \R\) in \eqref{eq:progdec} satisfy the following properties.
	\begin{enumeratass}
		\item \label{ass:progdec-global:3} 
		There exist parameters \(\mon, \com \in \R\)
		satisfying
		\(
			[\mon]_{-}[\com]_{-} < 1
		\)
		such that
		$S$ is $(\mon \proj_{X^\bot}, \com \proj_X)$\hyp{}semimonotone at 
		a solution $(x^\star, y^\star) \in \link_X S$.
		\item \label{ass:progdec-global:1} Operator \(S\) is outer semicontinuous. 
		\item \label{ass:progdec-global:2} The resolvent $J_{\gamma^{-1} S}$ has full domain.
	\end{enumeratass}
\end{assumption}

For instance, \Cref{ass:progdec-global} holds for any maximally $(\mon \proj_{X^\bot}, \com \proj_X)$\hyp{}semimonotone operator $S$ with nonempty solution set $\link_X S$.
Under this set of global conditions, \ref{eq:progdec} converges to a solution of linkage problem \eqref{prob:linkage-intro} from any initial point $(x^0, y^0) \in X \times X^\bot$, as stated below.

\begin{corollary}[global convergence of \ref{eq:progdec}]\label{thm:progdec}
    Suppose that \cref{ass:progdec-global} holds. Consider a sequence
	$(q^k, x^k, y^k)_{k\in\N}$
	generated by \ref{eq:progdec} starting from any $(x^0, y^0) \in X \times X^\bot$ with stepsize $\gamma$ and relaxation parameters $\lambda_x$ and $\lambda_y$ satisfying \eqref{eq:thm:progdec-local:stepsize}, and define
	with
	\(
		\bar x^k
			{}\coloneqq{}
		\proj_{X}(q^k)
	\)
	and
	\(
		\bar y^k
			{}\coloneqq{}
		y^k - \gamma\proj_{X^{\bot}}(q^k)
	\)
	as in \eqref{eq:progdec-full}.
    Then, either a point $(\bar{x}^k, \bar{y}^k) \in \link_X S$ is reached in a finite number of iterations or all the claims from \Cref{thm:progdec-local} hold.
\end{corollary}

\begin{table}
	\centering
	\caption{
		Existing and new convergence conditions for
		Spingarn's method of partial inverses
		(where $\gamma = \lambda_x = \lambda_y = 1$),  
		relaxed Douglas--Rachford splitting 
		(where $\lambda \coloneqq \lambda_x = \lambda_y$),
		standard progressive decoupling
		(where $\lambda_x = 1$, $\lambda_y = 1 - \nicefrac{[\mon]_-}{\gamma}$) 
		and progressive decoupling+ (see \ref{eq:progdec}) 
		applied to structured inclusion problem \eqref{prob:inclusion-intro} for the class of maximally $(\mon \proj_{X^\bot}, \com \proj_X)$\hyp{}semimonotone operators.
    }
    \label{tab:progdec:connections}
	\footnotesize
	\centering
	\ifsiam
\else
		\begin{tblr}{
				hlines,
				hline{3-5} = {1-3}{dotted},
				colspec={M{3.195cm}M{5.055cm}M{4.645cm}},
				cell{1}{1-3}={c,gray!20},
			}
			Method & Existing conditions & New conditions\\
			Spingarn's method
			& $\mon \geq 0$, $\com \geq 0$
			\cite[Thm. 4.1]{spingarn1983Partial}
			& $\mon > -\nicefrac12$, $\com > -\nicefrac12$ \\
			Relaxed DRS
			& 
			$\mon \geq 0$, $\com \geq 0$, $\gamma > 0$, $\lambda \in \bigl(0, 2\bigr)$ \cite[Thm. 7]{eckstein1992DouglasRachford}
			& 
			$[\mon]_-[\com]_- < 1$,
			$\gamma \in \bigl([\mon]_-$, $\nicefrac{1}{[\com]_-}\bigr)$, $\lambda \in \bigl(0, 2(1 + \min\{\gamma\com, \nicefrac\mon\gamma\})\bigr)$
			\\
			Progressive decoupling
			& $\mon \in \R$, $\com \geq 0$, $\gamma > [\mon]_-$
			\cite[Thm. 1]{rockafellar2019Progressive}
			& $[\mon]_-[\com]_- < \nicefrac12$, $\gamma \in \bigl([\mon]_-, \nicefrac{1}{2[\com]_-}\bigr)$\\
			Progressive decoupling+
			& 
			-
			& 
			$[\mon]_-[\com]_- < 1$,
			$\gamma \in \bigl([\mon]_-$, $\nicefrac{1}{[\com]_-}\bigr)$, $\lambda_x \in \bigl(0, 2(1 + \gamma\com)\bigr)$,
			$\lambda_y \in \bigl(0, 2(1 + \nicefrac\mon\gamma)\bigr)$
		\end{tblr}
	\fi
\end{table}

Recall that, by design, our algorithm \ref{eq:progdec} provides a unified framework that strictly generalizes Spingarn's method of partial inverses, progressive decoupling, and relaxed Douglas--Rachford splitting. Importantly, when we specialize our results to each of these methods, we not only recover but also strictly extend their known convergence guarantees. This is illustrated in \Cref{tab:progdec:connections}, which focuses on the maximally $(\mon \proj_{X^\bot}, \com \proj_X)$\hyp{}semimonotone case.
From this table, one can see that our new conditions on the moduli $\mon$ and $\com$, as well as on the algorithmic parameters, are significantly weaker than those required in prior work.
In particular, we relax the standard monotonicity assumption in Spingarn's method and relaxed DRS, and the elicitable monotonicity assumption in progressive decoupling.
Our results demonstrate that each of these methods remains applicable in the maximally $(\mon \proj_{X^\bot}, \com \proj_X)$\hyp{}semimonotone setting, albeit with potentially more limited nonmonotonicity---explicitly quantified by the moduli $\mon$ and $\com$.
As our analysis also covers the more general setting where semimonotonicity is neither maximal nor global,
we have thus developed a broad algorithmic framework for solving linkage problems of the form \eqref{prob:linkage-intro} under weak assumptions on the operator \(S\).

Finally, note that while much of the discussion in this work has centered on the nonmonotone regime---as it most clearly highlights the strength of our contributions---our framework also yields new ranges for the relaxation parameters $\lambda_x$ and $\lambda_y$ in the setting where $\mon$ and $\com$ are strictly positive, extending beyond the standard bound of two. 


  \subsection{Examples and sufficient conditions for continuously differentiable mappings}\label{subsec:progdec:examples}

We begin this subsection by presenting two specific applications of our convergence theory.
As demonstrated in \cite[Ex. 2.6]{evens2023convergence} through a simple linear example, the range of stepsize and relaxation parameters provided by \cref{thm:pppa} for \ref{eq:PPPA-intro} is tight.
Given the established equivalence between \ref{eq:progdec} and \ref{eq:PPPA-intro}, we expect the same tightness result to hold for the stepsize parameters of \ref{eq:progdec}.
This is illustrated in the following linear linkage problem, which 
can be interpreted as a worst-case example for the case where \(\lambda_x = \lambda_y\), i.e., the setting of relaxed DRS (recall \Cref{rem:equiv-progdec-DRS}).

\begin{example}[tightness of relaxation parameters]\label{ex:tightness}
    Consider linkage problem \eqref{prob:linkage-intro}, where 
    $
        X = \set{(x,0)}[x\in\R]
    $
    and
    $$
        S(x)
            {}\coloneqq{}
        \begin{bmatrix}
            1 + a^2 & 1\\
            1 & 1
        \end{bmatrix}\ \frac{1}{a} x
    $$
    for some $a \in \R\backslash\{0\}$.
    Consider a sequence
    $\seq{x^k, y^k}$
    generated by \ref{eq:progdec} (starting from $x^0 \in \R^2$, $y^0 \in \R^2$) with stepsize $\gamma = 1$ and relaxation parameter $\lambda = \lambda_x = \lambda_y > 0$.
    Then, the following hold.
    \begin{enumerate}
        \item\label{it:ex:tightness:1}
        \ifsiam\else
            The sequence
        \fi
        $\seq{x^k, y^k}$ converges to $\link_X S = (0,0)$ if and only if $\lambda$ lies in the interval
        \(
            (0, 2(1+\tfrac{a}{1 + a^2})).
        \)
        \item\label{it:ex:tightness:2}
        \ifsiam\else
            Operator
        \fi
        $S$ is $(\com \proj_{X^\bot}, \com \proj_X)$\hyp{}semimonotone where $\com = \tfrac{a}{1 + a^2}$.

        \item\label{it:ex:tightness:3}
        \cref{ass:progdec-global} holds and in \cref{ass:progdec-global:3} the set
		$ 
			\pazocal{S}^\star
		$
		is equal to
		$ 
			\link_X S
		$,
        and hence \Cref{thm:progdec} applies.
        Moreover, this result is tight in the sense that the bounds on $\lambda$ from condition \eqref{eq:thm:progdec-local:stepsize} for $\gamma = 1$ match the tight bounds from \cref{it:ex:tightness:1}.\qedhere
    \end{enumerate}
\end{example}

\begin{figure}
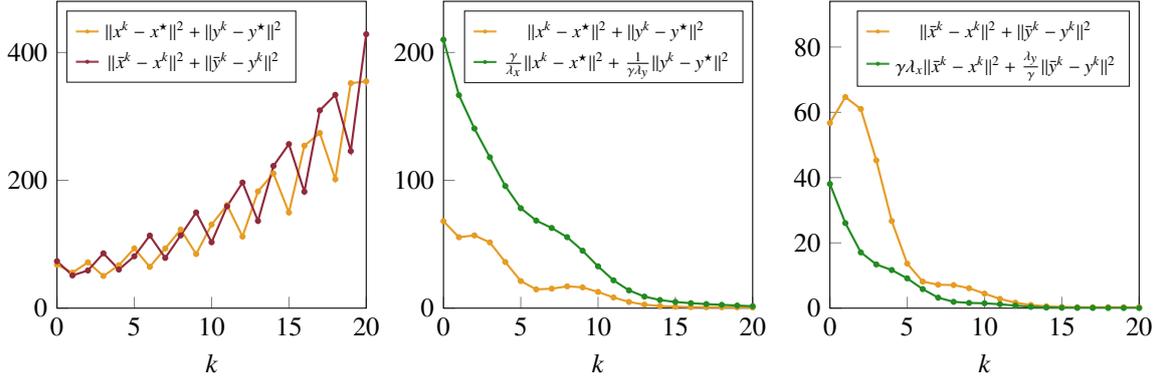

    \centering
    \ifsiam
\else
        \includetikz{Matrix-splitting/progdec-descent}
    \fi
    \caption{
        Different methods applied to \Cref{ex:linear-system}, starting from the pair $x_0 = (-2,-2,-2,-2) \in X$ and $y_0 = (1,1,-1,-1) \in X^\bot$.
        (left)
        \ref{eq:progdec} with $\gamma = \lambda_x = \lambda_y = 1$, i.e., 
        Spingarn's method of partial inverses.
        In this setting, both the sequence
        \(
            \seq{\nrm{x^k - x^\star}^2+ \nrm{y^k - y^\star}^2}
        \),
        representing the distance to the solution,
        and the sequence
        \(
            \seq{\nrm{\bar{x}^k - x^k}^2 + \nrm{\bar{y}^k - y^k}^2}
        \)
        do not converge to zero.
        (middle and right)
        \ref{eq:progdec} with $\gamma = \nicefrac{10}{9}$, $\lambda_x = \nicefrac95(1-\nicefrac\gamma2) = \nicefrac{4}{5}$ and $\lambda_y = \nicefrac95(1-\nicefrac1\gamma) = \nicefrac{9}{50}$, which complies with the stepsize conditions from \Cref{it:ex:linear-system:3}.
        As indicated by the convergence results from \Cref{it:progdec-local:v,it:progdec-local:rate}, the sequences
        \(
            \gamma\lambda_x^{-1}\nrm{x^k - x^\star}^2 + \gamma^{-1}\lambda_y^{-1}\nrm{y^k - y^\star}^2
        \)
        and
        \(
            \gamma\lambda_x\nrm{\bar{x}^k - x^k}^2+\gamma^{-1}\lambda_y\nrm{\bar{y}^k - y^k}^2
        \)
        are indeed nonincreasing and converge to zero.
        Note that although the sequences
        \(
            \seq{\nrm{x^k - x^\star}^2+ \nrm{y^k - y^\star}^2}
        \)
        and
        \(
            \seq{\nrm{\bar{x}^k - x^k}^2 + \nrm{\bar{y}^k - y^k}^2}
        \)
        also converge to zero, they are not nonincreasing.
    }
    \label{fig:nonsmooth:stationary}
\end{figure}

Our second example highlights the limitations of existing convergence theory.
Specifically, we consider solving a simple linear system of equations, formulated as a linkage problem over the consensus subspace. 
As shown in \cref{it:ex:linear-system:2}, the operator in this setup lacks elicitable monotonicity at any level, meaning it is not \((\mu \Pi_{X^\perp}, 0)\)-semimonotone for any \(\mu \in \mathbb{R}\).
Consequently, the convergence results for standard progressive decoupling, as discussed in \cite{rockafellar2019Progressive}, do not apply. However, we demonstrate that the operator is semimonotone and that our results ensure convergence of \ref{eq:progdec} under suitable stepsize conditions (see also \Cref{fig:nonsmooth:stationary}).

\begin{example}[linear system]\label{ex:linear-system}
    Consider solving the linear system
    \(
        b = M_1(x) + M_2(x),
    \)
    where
    \[
        b
            =
        \begin{bmatrix}
            2\\
            -3
        \end{bmatrix},
        \quad
        M_1
            =
        \begin{bmatrix}
            -1 & 2\\
            -2 & -1
        \end{bmatrix}
        \quad\text{and}\quad
        M_2
            =
        \begin{bmatrix}
            0 & 1\\
            0 & 0
        \end{bmatrix}.
    \]
    This is equivalent to finding an 
    $
        x \in X \coloneqq \set{(x_1,x_2) \in \R^2 \times \R^2}[x_1 = x_2],
    $
    such that
    \(b = M_1(x_1) + M_2(x_2)\).
    Using an auxiliary variable $y$ this can be cast as the linkage problem
    \begin{equation}
        \text{find} \quad (x, y) \in \R^4 \times \R^4 \quad \text{such that}
        \quad
        x
            \in
        X,
        \;
        y
            \in
        X^{\bot},
        \; 
        y
            = 
        S(x)
            \coloneqq
        \begin{bmatrix}
            M_1 & 0_{2\times 2}\\
            0_{2\times 2} & M_2
        \end{bmatrix}
        (x) - 
        \begin{bmatrix}
            0_2\\
            b
        \end{bmatrix}.
    \end{equation}
    The following hold.
    \begin{enumerate}
        \item\label{it:ex:linear-system:1} 
        \ifsiam\else
            Operator
        \fi
        $S$ is $(\mon \proj_{X^\bot}, \com \proj_X)$\hyp{}semimonotone, where $\mon = -1$ and $\com = \nicefrac{-1}2$.
        \item\label{it:ex:linear-system:2} 
        \ifsiam\else
            Operator
        \fi
        $S$ is not elicitable monotone, i.e., $S$ is not $(\mon \proj_{X^\bot}, 0)$\hyp{}semimonotone for any $\mon \in \R$.
        \item\label{it:ex:linear-system:3} 
        \cref{ass:progdec-global} holds and in \cref{ass:progdec-global:3} the set
		$ 
			\pazocal{S}^\star
		$
		is equal to
		$ 
			\link_X S
		$,
        and hence \Cref{thm:progdec} applies.
        In particular, any sequence $\seq{x^k, y^k}$ generated by \ref{eq:progdec} converges to
        $
            \link_X S
        $
        provided that
        \(
            \gamma \in (1,2)
        \),
        \(
            \lambda_x \in (0, 2(1-\nicefrac{\gamma}{2}))
        \)
        and
        \(
            \lambda_y \in (0, 2(1-\nicefrac1{\gamma}))
        \).\qedhere
    \end{enumerate}
\end{example}

 
In the remainder of this subsection, our goal is to gain more insights into the class of \((\mon \proj_{X^\bot}, \com \proj_X)\)\hyp{}semimonotone operators by providing sufficient conditions when dealing with linear operators, and more generally continuously differentiable operators.
We begin by stating an elementary fact about the semimonotonicity of linear mappings, which follows directly from the definition.

\begin{fact}[linear mappings]\label{fact:sm:linear}
    A linear mapping \(A \in \R^{n \times n}\) is $(\mon \proj_{X^\bot}, \com \proj_X)$\hyp{}semimonotone, 
    if and only if 
    \[
        \tfrac12(A + A^\top) -  \mu \proj_{X^\bot} - \rho A^\top \proj_{X} A \succeq 0,
    \]
\end{fact}


Remark that in the special case where \(\com = 0\), i.e., the elicitable monotone case studied in \cite{rockafellar2019Progressive}, this condition reduces to
\(
    \frac{1}{2}(A + A^\top) \succeq \mu \proj_{X^\bot},
\)
which may only hold when \(A\) is positive semidefinite on the linkage subspace \(X\), in the sense that
\(
    \proj_{X}(A + A^\top)\proj_{X} \succeq 0.
\)
In the setting considered in our work, where \(\rho \in \mathbb{R}\), \(A\) may be nonmonotone even on \(X\).


Next lemma establishes global semimonotonicity of a linear mapping \(A\), provided that condition \eqref{eq:lem:semimon:lin:cond} holds. Notably, the condition is required to hold only on the linkage subspace \(X\).
The proof relies primarily on \cref{fact:sm:linear} and the Schur complement lemma. 
In the framework of the linkage problem \cref{prob:linkage-intro} with operator \(S = A\), this result indicates that we only need to check condition \eqref{eq:lem:semimon:lin:cond} on the feasible subspace to establish global semimonotonicity.  

\begin{lemma}[global semimonotonicity from linearity] \label{lem:semimon:lin}
    Let \(X\subseteq \R^n\) be a closed linear subspace.
    Suppose that \(A\in\R^{n\times n}\) is a linear mapping for which there exist parameters $\mon > 0$ and $\com \in \R$ such that
    \begin{align*}
        \inner{x, Ax} \geq \mon \|x\|^2 + \com \|\proj_X Ax\|^2,
        \qquad \forall x \in X.
        \numberthis\label{eq:lem:semimon:lin:cond}
    \end{align*}
    Let 
    \begin{align*}
        \bar A
            \coloneqq
		\tfrac12(A + A^\top) - \com A^\top \proj_X A
        ,
        \quad
        \beta
          \dfn
        \nrm{\proj_{X}\bar A\proj_{X^\bot}}
        ,
        \quad
        \sigma
            {}\dfn{}
        \min_{y \in X^{\bot}, y \neq 0} \frac{\langle y, \bar A y \rangle}{\nrm{y}^2}.
        \numberthis\label{eq:lem:semimon:lin:params}
    \end{align*}
    Then, \(A\) is \(\Bigl(\bigl(\sigma-\nicefrac{\beta^2}{\mon}\bigr) \proj_{X^\bot}, \com \proj_X\Bigr)\)\hyp{}semimonotone, and maximally so if $\com \geq 0$.
    \begin{proof}
		By \Cref{fact:sm:linear}, \(A\) is \(\Bigl(\bigl(\sigma-\nicefrac{\beta^2}{\mon}\bigr) \proj_{X^\bot}, \com \proj_X\Bigr)\)\hyp{}semimonotone 
        if and only if
		\(
          \bar A \succeq \bigl(\sigma-\nicefrac{\beta^2}{\mon}\bigr) \proj_{X^\bot}.
        \) 
        Using the trivial identity \(v = \proj_X(v) + \proj_{X^\perp}(v)\) for any \(v\in \R^n\), it can be equivalently expressed as  
		\begin{align*}
            \begin{bmatrix}
                \proj_X \\
                \proj_{X^\bot}
            \end{bmatrix}^\top
            \begin{bmatrix}
                \proj_X \bar A \proj_X & \proj_X \bar A \proj_{X^\bot}\\
                \proj_{X^\bot} \bar A \proj_X & \proj_{X^\bot} \bar A \proj_{X^\bot}
                    {}-{}
                \bigl(\sigma-\nicefrac{\beta^2}{\mon}\bigr)\proj_{X^\bot}
            \end{bmatrix}
            \begin{bmatrix}
                \proj_X \\
                \proj_{X^\bot}
            \end{bmatrix}
                \succeq
            0.
            \numberthis\label{eq:lem:semimon:lin:basic}
        \end{align*}
        Since $\proj_X \bar A \proj_X \succeq \mon \proj_X \succeq 0$ owing to \eqref{eq:lem:semimon:lin:cond}, a sufficient condition for \eqref{eq:lem:semimon:lin:basic} is given by
        \begin{align*}
            \begin{bmatrix}
                \mon \proj_X & \proj_X \bar A \proj_{X^\bot}\\
                \proj_{X^\bot} \bar A \proj_X & \proj_{X^\bot} \bar A \proj_{X^\bot}
                    {}-{}
                \bigl(\sigma-\nicefrac{\beta^2}{\mon}\bigr)\proj_{X^\bot}
            \end{bmatrix}
                \succeq
            0.
            \numberthis\label{eq:lem:semimon:lin:suff}
        \end{align*}
        Note that $(\mon \proj_X)^\dagger = \nicefrac{1}{\mon} \proj_X$ and that
        \(
            \Bigl(I - (\mon \proj_X)(\mon \proj_X)^\dagger\Bigr)\proj_X \bar A \proj_{X^\bot}
                {}={}
            (I - \proj_X)\proj_X \bar A \proj_{X^\bot}
                {}={}
            0.
        \)
        Therefore, it follows from the Schur complement lemma \cite[Thm. 16.1]{gallier2011Geometric} that \eqref{eq:lem:semimon:lin:suff} holds if and only if
		\begin{align*}
			\proj_{X^\bot} \bar A \proj_{X^\bot}
				{}-{}
			\bigl(\sigma-\nicefrac{\beta^2}{\mon}\bigr)\proj_{X^\bot}
				{}-{}
			\nicefrac{1}{\mon} \proj_{X^\bot} \bar A \proj_X \bar A \proj_{X^\bot}
				{}\succeq{}
			0.
		\end{align*}
		Consequently, the claimed \(\bigl((\sigma-\nicefrac{\beta^2}{\mon}) \proj_{X^\bot}, \com \proj_X\bigr)\)\hyp{}semimonotonicity of \(A\) follows by definition of $\sigma$ and $\beta$, since $\proj_{X^\bot} \bar A \proj_{X^\bot} \succeq \sigma \proj_{X^\bot}$ and $\proj_{X^\bot} \bar A \proj_X \bar A \proj_{X^\bot} \preceq \beta^2 \proj_{X^\bot}$.
        Finally, the maximality claim for $\com \geq 0$ follows from \Cref{cor:structured-semimon:maximal}. \qedhere

    \end{proof}
\end{lemma}

In the special case where $\com = 0$, we return to the elicitable monotone setting from \cite{rockafellar2019Progressive},
and obtain the following corollary for the semimonotonicity of the sum of a linear mapping \(A\) and normal cone of a closed convex set \(C\).
This corollary relates to the setting where $T = F + N_C$ in \eqref{prob:linkage-intro}, which emerges for instance in the optimization setting when minimizing a function of the form \(f + \delta_C\) over the subspace \(X\).
Note that \Cref{cor:semimon:lin} provides a tighter \emph{elicitation level} 
\(
    \sigma - \nicefrac{\beta^2}{\mon}
\)
compared to \cite[Thm. 5]{rockafellar2019Progressive}, where the parameter \(\sigma\) is given by
\(
    \sigma = -\nrm{\proj_{X^\bot}A\proj_{X^\bot}}.
\)
In contrast, the parameter \(\sigma\) in \eqref{eq:cor:semimon:lin:params} is potentially larger, since
\[
    \sigma
        \geq
    \tfrac12\lambda_{\rm min}(\proj_{X^\bot}(A + A^\top)\proj_{X^\bot})
        \geq
    -\tfrac12\nrm{\proj_{X^\bot}(A + A^\top)\proj_{X^\bot}}
        \geq
    -\nrm{\proj_{X^\bot}A\proj_{X^\bot}}.
\]

\begin{corollary}\label{cor:semimon:lin}
    Let \(X\subseteq \R^n\) be a closed linear subspace and let $C \subseteq \R^n$ be a nonempty closed convex set.
    Suppose that \(A\in\R^{n\times n}\) is a linear mapping
    for which there exists a $\mon > 0$ such that
    \begin{align*}
        \inner{x, Ax} \geq \mon \|x\|^2,
        \qquad \forall x \in X.
    \end{align*}
    Define
    \begin{align*}
        \beta
            \dfn
        \tfrac12\nrm{\proj_{X}(A + A^\top)\proj_{X^\bot}}
        \quad\text{and}\quad
        \sigma
            \dfn
        \min_{y \in X^{\bot}, y \neq 0} \frac{\langle y, \proj_{X^\bot}(A + A^\top)\proj_{X^\bot} y \rangle}{2\nrm{y}^2}.
        \numberthis\label{eq:cor:semimon:lin:params}
    \end{align*}
    Then, \(T \coloneqq A + N_C\) is maximally \(\Bigl(\bigl(\sigma-\nicefrac{\beta^2}{\mon}\bigr) \proj_{X^\bot}, 0\Bigr)\)\hyp{}semimonotone.
    \begin{proof}
        The normal cone of a closed convex set is maximally monotone \cite[Ex. 20.26]{bauschke2017Convex}, and owing to \Cref{lem:semimon:lin} it holds that \(A - \bigl(\sigma-\nicefrac{\beta^2}{\mon}\bigr) \proj_{X^\bot}\) is maximally monotone.
        Consequently, it follows from \cite[Cor. 25.5(i)]{bauschke2017Convex} that \(T - \bigl(\sigma-\nicefrac{\beta^2}{\mon}\bigr) \proj_{X^\bot}\) is maximally monotone, establishing the claim.
    \end{proof}
\end{corollary}
Building on \Cref{lem:semimon:lin} for linear mappings, we now extend this result to general continuously differentiable mappings \(F : \R^n \rightarrow \R^n\).
The key idea is to express differences between evaluations of \(F\) as a weighted sum of Jacobian-vector products using the mean-value theorem \cite[Prop. 7.1.16]{facchinei2003FiniteDimensional} (see \eqref{eq:cor:semimon:lin:functions:global:mean-value}) and then apply \Cref{lem:semimon:lin} to the Jacobian matrices.
This leads to the following proposition, which provides sufficient conditions for a continuously differentiable mapping to be 
\(
    \bigl(
        \mon \proj_{X^\bot},
        \com \proj_X
    \bigr)
\)\hyp{}semimonotone. The full proof is deferred to \Cref{proof:cor:semimon:lin:functions}.

\begin{proposition}\label{cor:semimon:lin:functions}
    Let \(X\subseteq \R^n\) be a closed linear subspace and let $F : \R^n \rightarrow \R^n$ be a continuously differentiable mapping with Jacobians $\Jac F(x)$. Then, the following hold.
    \begin{enumerate}
        \item\label{cor:semimon:lin:functions:global-easy} 
        Let $\pazocal{U}^{\rm link}_1 \subseteq \R^n$ be a convex set and suppose that there exist scalars $\mon \in \R$ and $\com \in \R$ such that $\forall \tau \in (0,1)$
        \begin{align*}
            \inner{x - \other{x}, \Jac F\Bigl((1-\tau)x + \tau \other{x}\Bigr) (x - \other{x})} \geq \mon \|\proj_{X^\bot} (x - \other{x})\|^2 + \com \|\proj_X \Bigl(F(x) - F(\other{x})\Bigr)\|^2,
            \ifsiam\else
                \qquad \forall x, \other{x} \in \pazocal{U}^{\rm link}_1.
            \fi
        \end{align*} 
        Then, $F$ is
        \(
            \bigl(
                \mon \proj_{X^\bot},
                \com \proj_X
            \bigr)
        \)\hyp{}semimonotone on $\pazocal{U}^{\rm link}_1 \times \R^n$, and maximally so if $\pazocal{U}^{\rm link}_1 = \R^n$ and $\com \geq 0$.

        \item\label{cor:semimon:lin:functions:global} 
        Let $\pazocal{U}^{\rm link}_1 \subseteq \R^n$ be a convex set and suppose that there exists a function $\mon : \pazocal{U}^{\rm link}_1 \rightarrow \R_{++}$ satisfying
        \(
            \inf_{\tilde{x} \in \pazocal{U}^{\rm link}_1}
            \mon(\tilde{x}) > 0
        \)
        and a scalar $\com \geq 0$ such that
        \begin{align*}
            \inner{x, \Jac F(\tilde{x}) x} \geq \mon(\tilde{x}) \|x\|^2 + \com \|\proj_X \Jac F(\tilde{x})x\|^2,
            \qquad \forall x \in X, \forall \tilde{x} \in \pazocal{U}^{\rm link}_1.
            \numberthis\label{cor:semimon:lin:functions:global:ass}
        \end{align*}
        Then, defining the modulus
        \[
            \bar \mon \coloneqq \inf_{\tilde{x} \in \pazocal{U}^{\rm link}_1}
            \left\{
                \sigma(\tilde{x}) - \nicefrac{\beta^2(\tilde{x})}{\mon(\tilde{x})}
            \right\}
            \ifsiam\else
                \quad\text{where $\beta(\tilde{x})$ and $\sigma(\tilde{x})$ be defined as in \eqref{eq:lem:semimon:lin:params} with $\Jac F(\tilde{x})$ as matrix $A$},
            \fi
        \]
        it holds that $F$ is
        \(
            \bigl(
                \bar \mon \proj_{X^\bot},
                \com \proj_X
            \bigr)
        \)\hyp{}semimonotone on $\pazocal{U}^{\rm link}_1 \times \R^n$, and maximally so if $\pazocal{U}^{\rm link}_1 = \R^n$.

        \item\label{cor:semimon:lin:functions:local} Consider a point $\tilde{x} \in X$ 
        and suppose that there exist scalars $\mon > 0$ and $\com \geq 0$ such that
        \begin{align*}
            \inner{x, \Jac F(\tilde{x}) x} \geq \mon \|x\|^2 + \com \|\proj_X \Jac F(\tilde{x})x\|^2,
            \qquad \forall x \in X.
            \numberthis\label{cor:semimon:lin:functions:local:ass}
        \end{align*}
        Then, there exists a neighborhood $\pazocal{U}^{\rm link}_1$ of $\tilde{x}$ and scalars $\bar \mon \in \R, \bar \com \geq 0$ such that
        $F$ is
        \(
            \bigl(
                \bar \mon \proj_{X^\bot},
                \bar \com \proj_X
            \bigr)
        \)\hyp{}semimonotone on $\pazocal{U}^{\rm link}_1 \times \R^n$.
    \end{enumerate}
\end{proposition}
Observe that \cref{cor:semimon:lin:functions:global-easy} allows \(\com\) to be negative because the difference \( F(x) - F(\other{x}) \) remains inside the norm, avoiding the need to apply the mean-value theorem to this specific term.  
\cref{cor:semimon:lin:functions:global,cor:semimon:lin:functions:local} generalize \cite[Thm. 6]{rockafellar2019Progressive} to our more general setting with \(\com \geq 0\).  
Note that in the special case where \( F \) is linear, the first item of \Cref{cor:semimon:lin:functions} reduces to \Cref{fact:sm:linear}, while the second item then corresponds to \Cref{lem:semimon:lin} with \(\com \geq 0\).

An important special case of \Cref{cor:semimon:lin:functions} arises in the optimization setting (see \eqref{prob:linkage-optimization}), when \( F \) is the gradient of a twice continuously differentiable function \( f : \mathbb{R}^n \to \mathbb{R} \), meaning \( F = \nabla f \) with Hessian \( \Jac F = \nabla^2 f \).  
For instance, for any strict local minimum of the optimization problem \eqref{prob:linkage-optimization}, there exists a neighborhood \( \pazocal{U}^{\rm link}_1 \) around this minimum such that \( \nabla f \) is semimonotone on \( \pazocal{U}^{\rm link}_1 \times \mathbb{R}^n \), as summarized below.

\begin{example}[strict local minima]\label{ex:strict-local-min}
    Let \(X\subseteq \R^n\) be a closed linear subspace and let $U$ be an orthonormal basis for the range of $X$.
    Let $\psi : \R^n \rightarrow \R$ be a twice continuously differentiable function with gradients $\nabla \psi(x)$ and Hessians $\nabla^2 \psi(x)$.
    Let $x^\star \in X$ be a strict local minimum over the subspace \(X\), i.e., satisfying
    \(
        \nabla \psi(x^\star) \in X^\bot
    \)
    and $U^\top \nabla^2 \psi(x^\star) U \succ 0$.
    Then, condition \eqref{cor:semimon:lin:functions:local:ass} holds for some 
    $\mon > 0$ and $\com \geq 0$, 
    and by \cref{cor:semimon:lin:functions:local} there exists a neighborhood $\pazocal{U}^{\rm link}_1$ of $x^\star$ and a scalar $\bar \mon \in \R$ such that
    $\nabla \psi$ is
    \(
        \bigl(
            \bar \mon \proj_{X^\bot},
            \com \proj_X
        \bigr)
    \)\hyp{}semimonotone on $\pazocal{U}^{\rm link}_1 \times \R^n$.
    For example, one can take \(\mon = \lambda_{\rm min}\bigl(U^\top \nabla^2 \psi(x^\star) U\bigr)\) and \(\com = 0\).
\end{example}


We would like to highlight that similar results exist for nonsmooth functions, provided that they satisfy a \emph{variational convexity} \cite{rockafellar2019varconv} condition, which serves as a second-order sufficient condition for local optimality.
Specifically, if there exists a solution pair \((\tilde x, \tilde y) \in \link_X \partial f\) such that \(f - \frac{\mon}{2}\dist_X^2\) is variationally convex at \(\tilde x\) for \(\tilde y\), then \cite[Thm. 9]{rockafellar2019Progressive} guarantees that \(\tilde x\) is a local minimum over \(X\) and that a restriction of \(\gph \partial f\) is locally maximally \(\left(\mon \proj_{X^\bot}, 0\right)\)-semimonotone.
Extending this result to the case with nonzero \(\com\) remains an open question, which we leave for future work.

\section{Conclusion}\label{sec:conclusion}

In this work, we introduced a novel three-parameter algorithm tailored for solving nonmonotone linkage problems, which, for specific parameter choices, recovers Spingarn's method of partial inverses, progressive decoupling, and relaxed Douglas--Rachford splitting.
By interpreting our proximal algorithm as a halfspace projection, we established the local convergence of our algorithm and its specific instances, provided the involved operator belongs to a particular class of semimonotone operators.
Notably, the added parameter flexibility of our algorithm allows it to handle a higher degree of nonmonotonicity. 
To illustrate the tightness and practical significance of our results, we presented various exemplary problems throughout the paper.

Future research directions include exploring the degree of nonconvexity that our results allow within the framework of multistage stochastic programs compared to the standard progressive hedging algorithm, viewed through the lens of variational convexity.
Another main direction is to investigate the local convergence behavior of other splitting methods in the nonmonotone setting.

\begin{appendix}

\section{Auxiliary lemmas}\label{sec:auxiliary}
  
\begin{proposition}[partial inverse] 
	\label{prop:partialinverse:S-SX}
	Let \(\mon, \com \in \R\) and consider an operator \(S:\R^n\rightrightarrows\R^n\) and a closed linear subspace \(X\subseteq \R^n\).
	Then, the following are equivalent.
	\begin{enumerate}
		\item\label{it:prop:partialinverse:S} 
		\ifsiam\else
			Operator 
		\fi
		$S$ is (maximally) $\bigl(\mon\proj_{X^{\bot}}, \com\proj_{X}\bigr)$\hyp{}semimonotone \optional{at 
		$(x', y') \in \gph S$}
		on $\pazocal U$.
		\item\label{it:prop:partialinverse:SX}
		\ifsiam\else
			Operator 
		\fi
		\(S^X\) is (maximally) \(\bigl(\mon \proj_{X^\bot}+\com \proj_X \bigr)\)-comonotone \optional{at $L_X(x', y') \in \gph S^X$} on 
		\(
			L_X(\pazocal U).
		\)
	\end{enumerate}
	\begin{proof}
		Let
		$(x,y), (x', y') \in \gph S$
		and define
		\(
			(z, v)
				{}\coloneqq{}
			L_X(x, y)
				{}\in{}
			\gph S^X
		\)
		and
		\(
			(z', v')
				{}={}
			L_X(x', y')
				{}\in{}
			\gph S^X
		\).
		Then, it holds by definition of the partial inverse that 
		\[
			\langle
				z-z',v-v'
			\rangle
				{}={}
			\langle 
				x-x',y-y'
			\rangle
			\quad\text{and}\quad
			\qindef{v}{\mon \proj_{X^\bot}+\com \proj_X }
				{}={}
			\mon\|\proj_{X^{\bot}}(x - x')\|^{2}
				{}+{}
			\com\|\proj_{X}(y - y')\|^{2},
		\]
		establishing the claimed equivalence by definition of semimonotonicity.
	\end{proof}
\end{proposition}
  
\begin{corollary}\label{cor:structured-semimon:maximal}
    Let \(\mon \in \R\), \(\com \geq 0\) and consider a closed linear subspace \(X\subseteq \R^n\).
    If an operator \(S:\R^n\rightrightarrows\R^n\) is \((\mon\proj_{X^\bot}, \com\proj_{X})\)\hyp{}semimonotone and continuous, then it is also maximally \((\mon\proj_{X^\bot}, \com\proj_{X})\)\hyp{}semimonotone.

    \begin{proof}
        If \(S\) is \((\mon\proj_{X^\bot}, \com\proj_{X})\)\hyp{}semimonotone then it is also \(\bigl(\mon\proj_{X^\bot}, 0\bigr)\)\hyp{}semimonotone, as \(\com \geq 0\). By the continuity of \(S - \mon \proj_{X^\bot}\), it then follows from \cite[Cor. 20.28]{bauschke2017Convex} that \(S\) is maximally \(\bigl((\mon - \nicefrac{\com^2}{\mon}) \proj_{X^\bot}, 0\bigr)\)\hyp{}semimonotone.  
        The claim then follows from \cite[Prop. 4.2]{evens2023convergenceCP}. For completeness, we repeat the argument here.
        By maximality, there exists no \(\bigl((\mon - \nicefrac{\com^2}{\mon}) \proj_{X^\bot}, 0\bigr)\)\hyp{}semimonotone operator \(\tilde S\) such that \(\gph S \subset \gph \tilde S\).
        Since the class of \(\bigl((\mon - \nicefrac{\com^2}{\mon}) \proj_{X^\bot}, 0\bigr)\)\hyp{}semimonotone operators contains the class of \(\bigl((\mon - \nicefrac{\com^2}{\mon}) \proj_{X^\bot}, \com \proj_X\bigr)\)\hyp{}semimonotone operators for \(\com \geq 0\),
        the same conclusion holds for the latter class, completing the proof.
    \end{proof}
\end{corollary}

\section{Omitted proofs}\label{sec:omitted}
  \begin{appendixproof}{ex:Rosenbrock}[\ (Rosenbrock function)]
    The gradient 
    \(
        \nabla f(x,y) =
        (
            2x - 4bx(y-x^2),
            2b(y-x^2)
        )
    \)
    has a $(-\nicefrac1{16})$\hyp{}weak Minty solution at $(0,0)$ if and only if
    \(
        \inner*{
                (
                    x,
                    y
                ),
                \nabla f(x,y)
        }
            {}\geq{}
        -\tfrac1{16}                        
        \left\|\,             
            \nabla f(x,y)
        \right\|^2
    \)
    for all $x, y \in \R$,
    which is a quadratic inequality in $y$, equivalently given by
    \begin{align*}
        \left(b^2x^2+\tfrac14b^2 + 2b\right)y^2
            -
        \left(2b^2x^4 + \tfrac12b^2x^2 + 7bx^2\right)y
            +
        b^2x^6 + \tfrac14b^2x^4 + 5bx^4 + \tfrac94x^2 \geq 0.
    \end{align*}
    The coefficient for $y^2$ is strictly positive since $b > 0$. Therefore, the inequality holds for all $y \in \R$ since its discriminant is given by
    \(
        -bx^2\left(
            \tfrac94b + 18    
        \right)
    \),
    which is nonpositive for all $x \in \R$.
\end{appendixproof}
  \begin{appendixproof}{prop:local-comon:convex-solutions}
    Since the operator $\gph T \cap \pazocal U$ is $\DRSRho$\hyp{}comonotone, Zorn's lemma ensures that there exists a maximally \(\DRSRho\)\hyp{}comonotone operator $\hat T : \R^n \rightrightarrows \R^n$ satisfying
    \(
        \gph T \cap \pazocal U \subseteq \gph \hat T
    \).
    Moreover, since $T$ is maximally $\DRSRho$-cohypomonotone on $\pazocal U$ and $\hat T$ is $\DRSRho$-cohypomonotone, it follows that
    \(
        \gph T \cap \pazocal U = \gph \hat T \cap \pazocal U
    \), 
    which directly implies that
    \(
        \zer T \cap \pazocal U = \zer \hat T \cap \pazocal U
    \).
    Therefore, it suffices to proof that \(\zer \hat T\) is convex, since the intersection of two convex sets is convex.

    Observe that
    $
        \zer\, (\hat T^{-1}+\DRSRho)^{-1}
            =
        (\hat T^{-1}+\DRSRho)(0)
            = 
        \zer \hat T
    $.
    Furthermore, $T^{-1}+\DRSRho$ is maximally monotone by definition, and thus so is $(T^{-1}+\DRSRho)^{-1}$.
    Since the set of zeros of a maximally monotone operator is convex \cite[Prop. 23.39]{bauschke2017Convex}, \(\zer \hat T\) is convex and the proof is completed.
\end{appendixproof}

\begin{appendixproof}{prop:local-comon}
    \Cref{ass:PPPA-local:1} holds by maximal 
    \(\DRSRho\)\hyp{}comonotonicity of \(T\) on $\pazocal U_1 \times \pazocal U_2$.
    Let $(x^k, y^k) \in \graph T \cap (\pazocal U_1 \times \pazocal U_2)$ where \(x^k \rightarrow \bar{x}\) and $y^k \rightarrow \bar{y}$, and note that $x^k \in (T^{-1} - \DRSRho)^{-1}(y^k - \DRSRho x^k)$.
    Since $T$ is maximally $\DRSRho$\hyp{}comonotone on $\pazocal U_1 \times \pazocal U_2$, $(T^{-1} - \DRSRho)^{-1}$ is maximally monotone on $\pazocal U_1 \times \pazocal U_2$ and thus also outer semicontinuous \cite[Ex. 12.8]{rockafellar2009Variational}. Consequently, it follows that $\bar{x} \in (T^{-1} - \DRSRho)^{-1}(\bar{y} - \DRSRho \bar{x})$.
    This implies that \(\bar{y} \in T(\bar{x})\), showing outer semicontinuity of \(T\), i.e., that \Cref{ass:PPPA-local:0} holds.

    By construction, the operator $\gph T \cap (\pazocal U_1 \times \pazocal U_2)$ is $\DRSRho$\hyp{}comonotone.
    Therefore, owing to Zorn's lemma, there exists a maximally \(\DRSRho\)\hyp{}comonotone operator $\hat T : \R^n \rightrightarrows \R^n$ such that
    \begin{align*}
        \gph T \cap \left(\pazocal U_1 \times \pazocal U_2\right) \subseteq \gph \hat T.
        \numberthis\label{eq:prop:local-comon:gph}
    \end{align*}
    Then, owing to \Cref{ass:PPPA-local:2} the preconditioned resolvent 
    \((\M + \hat T)^{-1} \circ \M\) has full domain \cite[Prop. A.2]{evens2023convergence}.
    Let $z \in \pazocal U_1$ and consider the points
    \(
        \bar z \in (\M + \hat T)^{-1} \M(z)
    \)
    and 
    \(
        z^+ \coloneqq z + \Lambda (\bar z - z).
    \)
    Since $(\bar z, \M(z - \bar z)) \in \gph \hat T$ it holds by $\DRSRho$\hyp{}comonotonicity of $\hat T$
    for any $z^\star \in \pazocal{S}^\star \subseteq \zer \hat T$ that
    \ifsiam
\else
        \begin{align*}
            \inner{\M(\bar z - z), z - z^\star}
                {}={}
            - \|\bar z - z\|_{\M}^2 - \inner{P(z-\bar z), \bar z - z^\star}
                {}\leq{}
            - \|\bar z - z\|_{\M^{\nicefrac12}(\I + \M^{\nicefrac12}\DRSRho\M^{\nicefrac12})\M^{\nicefrac12}}^2
                {}\leq{}
            -\bar \alpha\|\bar z - z\|_{\M \Lambda},  
            \numberthis\label{eq:prop:local-comon:inner}
        \end{align*}
    \fi
    where \(\bar \alpha > \nicefrac12\) is defined as in \eqref{eq:half:alpha}.
    This implies that
    \begin{align*}
        \|z^+ - z^\star\|_{\M \Lambda^{-1}}^2
            {}={}
        \|z - z^\star + \Lambda (\bar z - z)\|_{\M \Lambda^{-1}}^2
            {}={}&
        \|z - z^\star\|_{\M \Lambda^{-1}}^2 + 2\inner{\bar z - z, z - z^\star}_{\M} + \|\bar z - z\|_{\M \Lambda}^2\\
            \dueto{\eqref{eq:prop:local-comon:inner}}
            {}\leq{}&
        \|z - z^\star\|_{\M \Lambda^{-1}}^2
            -
        (2\bar\alpha-1)\|\bar z - z\|_{\M \Lambda}^2\\
            {}\leq{}&
        \|z - z^\star\|_{\M \Lambda^{-1}}^2.
        \numberthis\label{eq:prop:local-comon:zplus}
    \end{align*}
    Consequently, using \(\bar z = z + \Lambda^{-1}(z^+ - z)\) and the triangle inequality, we obtain that
    \ifsiam
\else
        \begin{align*}
            \dist_{\M \Lambda^{-1}}(\bar z, \pazocal{S}^\star)
                {}\leq{}
            \|\bar z - z^\star\|_{\M \Lambda^{-1}}
                {}={}
            \|z - z^\star + \Lambda^{-1}(z^+ - z)\|_{\M \Lambda^{-1}}
                {}={}&
            \|(\I - \Lambda^{-1})(z - z^\star) + \Lambda^{-1}(z^+ - z^\star)\|_{\M \Lambda^{-1}}\\
                {}\leq{}&
            \|\I - \Lambda^{-1}\|\,\|z - z^\star\|_{\M \Lambda^{-1}} + \|\Lambda^{-1}\|\,\|z^+ - z^\star\|_{\M \Lambda^{-1}}\\
                \dueto{\eqref{eq:prop:local-comon:zplus}}{}\leq{}&
            \Bigl(\|\I - \Lambda^{-1}\| + \|\Lambda^{-1}\|\Bigr)\|z - z^\star\|_{\M \Lambda^{-1}}
            \numberthis\label{eq:prop:local-comon:barz}
        \end{align*}
    \fi
    and
    \ifsiam
\else
        \begin{align*}
            \dist_{\M \Lambda^{-1}}(\M(z - \bar z), 0)
                {}={}
            \|\M(z - \bar z)\|_{\M \Lambda^{-1}}
                {}\leq{}
            \|\M\| \|z - \bar z\|_{\M \Lambda^{-1}}
                {}\leq{}&
            \|\M\|\left(
                \|z - z^\star\|_{\M \Lambda^{-1}}  
                    +  
                \|\bar z - z^\star\|_{\M \Lambda^{-1}}  
            \right)\\
                \dueto{\eqref{eq:prop:local-comon:barz}}
                {}\leq{}&
            \|\M\|\left(
                1 + \|\I - \Lambda^{-1}\| + \|\Lambda^{-1}\|
            \right)
            \|z - z^\star\|_{\M \Lambda^{-1}}.
            \numberthis\label{eq:prop:local-comon:v}
        \end{align*}
    \fi
    Let $z^\star = \proj_{\pazocal{S}^\star}^{\M \Lambda^{-1}}(z)$ in \eqref{eq:prop:local-comon:barz} and \eqref{eq:prop:local-comon:v}, so that $\|z - z^\star\|_{\M \Lambda^{-1}} = \dist_{\M \Lambda^{-1}}(z, \pazocal{S}^\star)$.
    Then, for any 
    \[
        z
            {}\in{}
        \pazocal W \coloneqq \set{z \in \R^n}[\dist_{\M \Lambda^{-1}}(z, \pazocal{S}^\star) \leq         
        \min
        \left(
            \ifsiam
\else
                \frac{\delta_1}{\|\I - \Lambda^{-1}\| + \|\Lambda^{-1}\|},
                \frac{\delta_2}{
                \|\M\|\left(
                    1 + \|\I - \Lambda^{-1}\| + \|\Lambda^{-1}\|
                \right)}
            \fi
        \right)],
    \]
    it holds by 
    \eqref{eq:prop:local-comon:barz} that 
    $\bar z \in \pazocal U_1$ and by \eqref{eq:prop:local-comon:v} that $\M(z - \bar z) \in \pazocal U_2$.
    Since $T$ is maximally $\DRSRho$-cohypomonotone on $\pazocal U_1 \times \pazocal U_2$ and since $\hat T$ is $\DRSRho$-cohypomonotone, it follows from \eqref{eq:prop:local-comon:gph} that 
    \(
        \gph T \cap (\pazocal U_1 \times \pazocal U_2) = \gph \hat T \cap (\pazocal U_1 \times \pazocal U_2)
    \).
    Consequently, it holds that
    \(
        (\bar z, P(z - \bar z)) \in \gph T \cap (\pazocal U_1 \times \pazocal U_2)
    \),
    showing \Cref{ass:PPPA-local:0.5}.

    Finally, when $\pazocal U_1 \times \pazocal U_2 = \R^n \times \R^n$, the argument for \Cref{ass:PPPA-global:1,ass:PPPA-global:0}
    is analogous to the one above for
    \Cref{ass:PPPA-local:1,ass:PPPA-local:0}.
    To establish that the preconditioned resolvent 
    \((\M + T)^{-1} \circ \M\) has full domain, i.e., \cref{ass:PPPA-global:0.5}, we can directly apply \cite[Prop. A.2]{evens2023convergence}, completing the proof.
\end{appendixproof}
  \begin{appendixproof}{ex:Rosenbrock:progdec}[\ (Rosenbrock function)]
    Note that 
    \begin{align*}
        \nabla f(x)
            {}={}
        \begin{bmatrix}
            x_2 - 4bx_1(x_3 - x_1^2)\\
            x_1\\
            2b(x_3 - x_1^2)
        \end{bmatrix},
        \quad
        \proj_X
            {}={}
        \begin{bmatrix}
            \nicefrac12 & \nicefrac12 & 0\\
            \nicefrac12 & \nicefrac12 & 0\\
            0 & 0 & 1
        \end{bmatrix}
        \quad\text{and}\quad
        \proj_{X^\bot}
            {}={}
        \begin{bmatrix}
            \nicefrac12 & -\nicefrac12 & 0\\
            -\nicefrac12 & \nicefrac12 & 0\\
            0 & 0 & 0
        \end{bmatrix}.
    \end{align*}
    Therefore, $\nabla f$ is 
    $\Bigl(-\nicefrac94 \proj_{X^\bot}, -\nicefrac{1}{4}\proj_X\Bigr)$-semimonotone
    at $(x^\star, 0)$ if and only if
    \begin{align*}
        \inner{x, \nabla f(x)} \geq -\nicefrac94\|\proj_{X^\bot}x\|^2 - \nicefrac{1}{4} \|\proj_X \nabla f(x)\|^2
        \ifsiam
\else
            ,
        \fi
    \end{align*}
    \ifsiam
\else
        which is equivalent to a quadratic inequality in $x_2$, given by
        \begin{align*}
            \tfrac54x_2^2 
                -
            b\left(
                x_3-x_1^2
            \right)x_1x_2
                +
            2b^2x_1^6
                +
            \left(
                    -
                4b^2x_3
                    +
                b^2 
                    +
                5b
            \right)
            x_1^4
                +
            \left(
                2b^2x_3^2 
                    -
                2b^2x_3
                    -
                7bx_3 
                + 
                \tfrac54
            \right)x_1^2
                +
            \left(b^2+2b\right)x_3^2
                \geq
            0.
        \end{align*}
        This inequality holds for all $x_2 \in \R$ if and only if its discriminant is nonpositive, i.e., if and only if
        \begin{align*}
            b^2(
                x_3-x_1^2
            )^2x_1^2
                -
            5\left(
                2b^2x_1^6
                    +
                \left(
                        -
                    4b^2x_3
                        +
                    b^2 
                        +
                    5b
                \right)
                x_1^4
                    +
                \left(
                    2b^2x_3^2 
                        -
                    2b^2x_3
                        -
                    7bx_3 
                    + 
                    \tfrac54
                \right)x_1^2
                    +
                \left(b^2+2b\right)x_3^2
            \right)
                \leq
            0.
        \end{align*}
        By reordering the terms, this discriminant condition reduces to another quadratic inequality in $x_3$, given by
        \begin{align*}
            \left(
                -9b^2x_1^2
                    -
                5b^2
                    -
                10b
            \right)
            x_3^2
                +
            \left(
                18b^2x_1^4
                    +
                10b^2x_1^2
                    +
                35bx_1^2
            \right)x_3
                -
            9b^2x_1^6
                -
            (5b^2
                +
            25b)x_1^4
                -
            \tfrac{25}{4}x_1^2
                \leq
            0.
        \end{align*}
        This inequality holds for all $x_3 \in \R$ since its discriminant is given by
        \(
            -125b x_1^2 (2+b)
        \),
        which is nonpositive for all $x_1 \in \R$, completing the proof.
    \fi
\end{appendixproof}
  \begin{appendixproof}{ex:double-well}[\ (local minimum)]
    The claim is equivalent to showing for all $x = (x_1, x_2) \in \pazocal U$ that
    \(
        \inner{x, \nabla f(x)} \geq \|\proj_X \nabla f(x)\|^2.
    \)
    Plugging in the definition of $f$ and reordering, this is equivalent to
    \(
        x_1^2 + x_2^2 - x_1^2 x_2^2 \geq \tfrac18(x_1+x_2)^2(2-x_1x_2)^2,
    \)
    which can be shown to hold for all $(x_1, x_2) \in \R^2$ satisfying $x_1^2 + x_2^2 \leq 4$, completing the proof.
\end{appendixproof}
  \begin{appendixproof}{ex:consensus:cond}[\ (matrix splitting)]
    By assumption, it holds that  
    \(
        \tfrac12(A_i + A_i^\top) \succeq \mu \I + \rho A_i^\top A_i.
    \)
    Due to the block diagonal structure of \( A \) and condition \eqref{eq:ex:linear:StS0}, this implies that  
    \begin{align*}
        \tfrac12(A + A^\top)
            \succeq
        \mu \I + \rho A^\top A
            \overrel[\succeq]{\eqref{eq:ex:linear:StS0}}
        \mu \I + \frac{\rho}{\nu} A^\top 1 1^\top A.
        \numberthis\label{eq:ex:consensus:cond:given}
    \end{align*}
    For the consensus constraint, it holds that  
    \(
        \proj_{X} = \tfrac1N 1 1^\top
    \). Substituting this into \eqref{eq:ex:consensus:cond:given} gives  
    \begin{align*}
        \tfrac12(A + A^\top)
            \succeq
        \mon \proj_{X^\bot} + \mon \proj_{X} + \frac{N\rho}{\nu} A^\top \proj_{X} A.
    \end{align*}
    The claim then follows directly from \Cref{fact:sm:linear}, using that \(\mon \proj_{X} \succeq 0\).
\end{appendixproof}
  \begin{appendixproof}{prop:local-semi}
    \Cref{ass:progdec-local:3} holds since \(S\) is maximally 
    $(\mon \proj_{X^\bot}, \com \proj_X)$\hyp{}semimonotone on
    \(
        \pazocal{U}^{\rm link}
    \)
    and \(\link_X S\) is assumed to be nonempty.
    Let $(x^k, y^k) \in \graph S \cap (\pazocal{U}^{\rm link}_1 \times \pazocal{U}^{\rm link}_2)$
    where
    \(x^k \rightarrow \bar{x}\)
    and
    \(y^k \rightarrow \bar{y}\).
    By construction, this is equivalent to
    \(
        x^k - \com \proj_X y^k
            \in
        \bar{S}(y^k - \mon \proj_{X^\bot} x^k)
    \),
    where
    \(
        \bar{S}
            \coloneqq
        \left((S - \mon \proj_{X^\bot})^{-1} - \com \proj_X\right)^{-1}.
    \)
    Since $S$ is maximally 	
    $(\mon \proj_{X^\bot}, \com \proj_X)$\hyp{}semimonotone on
    \(
        \pazocal{U}^{\rm link}
    \), 
    \(
        \bar{S}
    \)
    is maximally monotone on
    \(
        \pazocal{U}^{\rm link}
    \)
    and thus also outer semicontinuous \cite[Ex. 12.8]{rockafellar2009Variational}. Consequently, it follows that
    \(
        \bar{x} - \com \proj_X \bar{y}
            \in
        \bar{S}(\bar{y} - \mon \proj_{X^\bot} \bar{x}).
    \)
    This implies that \(\bar{y} \in S(\bar{x})\), showing outer semicontinuity of \(S\), i.e., that \Cref{ass:progdec-local:1} holds.
    Since it is assumed that $\gamma$, $\lambda_x$ and $\lambda_y$ satisfy
    \eqref{eq:thm:progdec-local:stepsize}, it only remains to show that there exists an $\varepsilon > 0$ such that \Cref{ass:progdec-local:2} holds.
    To this end, define the sets
    \[
        \pazocal X_1
            {}\coloneqq{}
        \set{x \in X}[\|x - x^\star\|\leq \frac{\delta}{\sqrt{2}}]
        \quad\text{and}\quad
        \pazocal X_2
            {}\coloneqq{}
        \set{y \in X^\bot}[\|y\| \leq \frac{\delta}{\sqrt{2}}].
    \]
    Then, by construction
    \(
        \pazocal X
            {}\coloneqq{}
        \pazocal X_1+ \pazocal X_2	
            {}\subset{}
        \pazocal{U}^{\rm link}_1
            {}={}
        \set{x + y}[x \in X, y \in X^\bot, \|x- x^\star\|^2 + \|y\|^2 \leq \delta^2]
    \),
    since $x^\star \in X$.
    Analogously, defining
    \[
        \pazocal Y_1
            {}\coloneqq{}
        \set{x \in X}[\|x\| \leq \frac{\delta}{\sqrt{2}}]
        \quad\text{and}\quad
        \pazocal Y_2
            {}\coloneqq{}
        \set{y \in X^\bot}[\|y - y^\star\| \leq \frac{\delta}{\sqrt{2}}],
    \]
    it follows that
    \(
        \pazocal Y
            {}\coloneqq{}
        \pazocal Y_1 + \pazocal Y_2
            {}\subset{}
        \pazocal{U}^{\rm link}_2.
    \)
    Therefore, it follows that \(S\) is maximally $(\mon \proj_{X^\bot}, \com \proj_X)$\hyp{}semimonotone on
    \(
        \pazocal X \times \pazocal Y
            {}\subset{}
        \pazocal{U}^{\rm link}_1 \times \pazocal{U}^{\rm link}_2
    \).
    Applying \cref{prop:partialinverse:S-SX}, this is equivalent to $S^X$ being maximally \(\bigl(\mon \proj_{X^\bot}+\com \proj_X \bigr)\)\hyp{}comonotone on
    \(
        L_X(\pazocal X, \pazocal Y)
    \).
    By definition of the partial inverse and the sets \(\pazocal X\) and \(\pazocal Y\) it holds that
    \(
        L_X(\pazocal X, \pazocal Y)
            =
        (\pazocal X_1+ \pazocal Y_2) \times (\pazocal X_2 + \pazocal Y_1)
    \).
    Define the point $z^\star \coloneqq x^\star + y^\star$ and the matrices
    \(
        \M \coloneqq \gamma \proj_X + \gamma^{-1} \proj_{X^\bot}
    \)
    and
    \(
        \Lambda \coloneqq \lambda_x \proj_X + \lambda_y \proj_{X^\bot}
    \),
    for which
    \(
        \lambda_{\rm min}\bigl(\M \Lambda^{-1}\bigr)
            {}={}
        \min\bigl\{\gamma\lambda_x^{-1}, (\gamma\lambda_y)^{-1}\bigr\} > 0
    \).
    Then, it holds that
    \ifsiam
\else
        \begin{align*}
            \pazocal U_1
                {}\coloneqq{}&
            \set{z \in \R^n}[\|z - z^\star\|_{\M \Lambda^{-1}} \leq \frac{\delta}{\sqrt{2}}\sqrt{\lambda_{\rm min}\bigl(\M \Lambda^{-1}\bigr)}]\\
                {}\subseteq{}&
            \set{z \in \R^n}[\|z - z^\star\| \leq \frac{\delta}{\sqrt{2}}]
                {}={}
            \set{x + y}[x \in X, y \in X^\bot, \|x - x^\star\| + \|y-y^\star\| \leq \frac{\delta}{\sqrt{2}}]
                {}\subset{}
            \pazocal X_1+ \pazocal Y_2,
        \end{align*}
    \fi
    and analogously
    \begin{align*}
        \pazocal U_2
            {}\coloneqq{}&
        \set{v \in \R^n}[\|v\|_{\M \Lambda^{-1}} \leq \frac{\delta}{\sqrt{2}}\sqrt{\lambda_{\rm min}\bigl(\M \Lambda^{-1}\bigr)}]\\
            {}\subseteq{}&
        \set{v \in \R^n}[\|v\| \leq \frac{\delta}{\sqrt{2}}]
            {}={}
        \set{x + y}[x \in X, y \in X^\bot, \|x\| + \|y\| \leq \frac{\delta}{\sqrt{2}}]
            {}\subset{}
        \pazocal X_2 + \pazocal Y_1.
    \end{align*}
    Since we already established that $S^X$ is maximally \(\bigl(\mon \proj_{X^\bot}+\com \proj_X \bigr)\)\hyp{}comonotone on
    \(
        (\pazocal X_1+ \pazocal Y_2) \times (\pazocal X_2 + \pazocal Y_1)
    \),
    it follows that $S^X$ is also maximally \(\bigl(\mon \proj_{X^\bot}+\com \proj_X \bigr)\)\hyp{}comonotone on $\pazocal U_1 \times \pazocal U_2$ and we can invoke \cref{prop:local-comon} with $\delta_1 = \delta_2 = \nicefrac{\delta\sqrt{\lambda_{\rm min}\bigl(\M \Lambda^{-1}\bigr)}}{\sqrt{2}}$
    to show that \cref{ass:PPPA-local} holds with
    \(
        \pazocal U
            =
        \pazocal U_1 \times \pazocal U_2
    \)
    and
    \[
        \varepsilon
            {}={}
        \frac{
            \delta\sqrt{\lambda_{\rm min}\bigl(\M \Lambda^{-1}\bigr)}			
            \min
            \left(
                1,
                \|\M\|^{-1}
            \right)
        }{
            \sqrt{2}
            \left(
                \|\I - \Lambda^{-1}\| + \|\Lambda^{-1}\|
            \right)
        }
            {}={}
        \frac{
            \delta\sqrt{\min\bigl\{\gamma\lambda_x^{-1}, (\gamma\lambda_y)^{-1}\bigr\}}			
            \min(\gamma, \gamma^{-1})
        }{
            \sqrt{2}\left(
                1 - \min\{\lambda_x^{-1}, \lambda_y^{-1}\} + \max\{\lambda_x^{-1}, \lambda_y^{-1}\}
            \right)
        }.
    \]
    Therefore, by \Cref{lem:equiv-ass}, \Cref{ass:progdec-local:2} holds with the same \( \varepsilon \) and  
    \(
        \pazocal{U}^{\rm link} = L_X(\pazocal{U}_1 \times \pazocal{U}_2).
    \)  
    Since we already established that
    \(
        L_X(\pazocal{U}_1 \times \pazocal{U}_2) \subset \pazocal{U}^{\rm link}_1 \times \pazocal{U}^{\rm link}_2,
    \) 
    \Cref{ass:progdec-local:2} also holds for the larger set
    \(
        \pazocal{U}^{\rm link} = \pazocal{U}^{\rm link}_1 \times \pazocal{U}^{\rm link}_2,
    \)  
    which completes the first part of the proof.

    When $\pazocal{U}^{\rm link}_1 \times \pazocal{U}^{\rm link}_2 = \R^n \times \R^n$, then it follows from \cref{prop:partialinverse:S-SX} that $S^X$ is maximally \(\bigl(\mon \proj_{X^\bot}+\com \proj_X \bigr)\)\hyp{}comonotone on $\pazocal U_1 \times \pazocal U_2 = \R^n \times \R^n$.
    Therefore, \Cref{ass:PPPA-local} holds globally owing to \cref{prop:local-comon},
    and the claim follows by taking $\pazocal U^{\rm link} = \R^n \times \R^n$ in \cref{lem:equiv-ass}.
\end{appendixproof}
  \begin{appendixproof}{ex:tightness}[\ (tightness of relaxation parameters)]
    \begin{proofitemize}
        \item\ref{it:ex:tightness:1}:
        Note that $\proj_X = \diag(1,0)$ and $\proj_{X^\bot} = \diag(0,1)$, so that by \cref{it:spingarn:linear} 
        the partial inverse of $S$ with respect to $X$ is given by
        \begin{align*}
            S^X
                {}={}
            (\proj_{X^\bot} + \proj_X S)(\proj_X + \proj_{X^\bot} S)^{-1}
                {}={}&
            \left(
            \begin{bmatrix}
                1 + a^2 & 1\\
                0 & a
            \end{bmatrix}\
            \frac{1}{a}
            \right)\left(
                \begin{bmatrix}
                    a & 0\\
                    1 & 1
                \end{bmatrix}\
                \frac{1}{a}
            \right)^{-1}
                {}={}
            \begin{bmatrix}
                a & 1\\
                -1 & a
            \end{bmatrix}.
        \end{align*}
        Since $\zer S^X = \set{0}$, this implies by \cref{it:spingarn:zero} that $\link_X S = (0,0)$. Define
        \[
            H
                {}\coloneqq{}
            \I_2 + \lambda\bigl((\I + S^X)^{-1} - \I_2\bigr)
                {}={}
            \I_2 + \lambda\left(  
            \frac{1}{(1+a)^2+1}          
                \begin{bmatrix}
                1+a & -1\\
                1 & 1+a
            \end{bmatrix} - \I_2\right).
        \]
        By \Cref{lem:equivalence:PPPA-progdec} each iteration of \ref{eq:progdec} can be expressed as the linear dynamical system
        \(
            z^{k+1} = Hz^k
        \),
        where 
        \(
            x^k = \proj_X(z^k)
        \)
        and
        \(
            y^k = \proj_{X^\bot}(z^k)
        \).
        This linear system is globally asymptotically stable if and only if 
        the spectral radius of $H$, given by
        \(
            \sqrt{1-\tfrac{\lambda (2(1 + a + a^2) - \lambda(1 + a^2))}{(1 + a)^2 + 1}}
        \),
        is strictly less than one, which holds if and only if $\lambda$ lies in the interval
        \(
            (0, 2(1+\tfrac{a}{1 + a^2})).
        \)
        \item \ref{it:ex:tightness:2}:
        \ifsiam
\else
            Holds by definition of semimonotonicity,
        \fi
        since
        \(
            \tfrac12(S + S^\top)
                =
            \com \left(\proj_{X^\bot} + S^\top \proj_X S\right)
                =
            \begin{bsmallmatrix}
                1 + a^2 & 1\\
                1 & 1
            \end{bsmallmatrix}\ \frac1a.
        \)
        \item \ref{it:ex:tightness:3}:
        Note that $[\rho]_-[\rho]_- = \rho^2 = \frac{a^2}{1+a^2} < 1$. Therefore, it follows directly from \cref{it:ex:tightness:2} that \cref{ass:progdec-global:3} holds.
        Moreover, operator $S$ is continuous and thus outer semicontinuous, showing \cref{ass:progdec-global:1}.
        Finally, note that $J_{\gamma^{-1} S} = (\I + S)^2$ is invertible and hence the resolvent has full domain, showing \cref{ass:progdec-global:2}, and is single-valued. \qedhere
    \end{proofitemize}
\end{appendixproof}
  \begin{appendixproof}{ex:linear-system}[\ (linear system)]
    \begin{proofitemize}
        \item \ref{it:ex:linear-system:1}:
        Let 
        \(
            M
                {}\coloneqq{}
            \begin{bsmallmatrix}
                M_1 & 0_{2\times 2}\\
                0_{2\times 2} & M_2
            \end{bsmallmatrix}
        \)
        and
        \(
            m
                {}\coloneqq{}
            \begin{bsmallmatrix}
                0_2\\
                b
            \end{bsmallmatrix}
        \),
        so that $S(x) \coloneqq Mx - m$.
        Then, by the definition of semimonotonicity, the claim for $S$ holds if and only if
        \[  
            \bar M
                {}\coloneqq{}
            \tfrac12(M + M^\top)
                -
            \mon \proj_{X^\bot} - \com M^\top \proj_X M
                =
            \frac14
            \begin{bmatrix}
                3 & 0 & -2 & -1\\
                0 & 3 & 0 & 0\\
                -2 & 0 & 2 & 2\\
                -1 & 0 & 2 & 3
            \end{bmatrix}
                \succeq
            0.
        \]
        where we used that the projections onto $X$ and $X^\bot$ are given by
        \(
            \proj_X 
                {}={}
            \frac12
            \begin{bsmallmatrix}
                \I_2 & \I_2\\
                \I_2 & \I_2
            \end{bsmallmatrix}
        \)
        and
        \(
            \proj_{X^\bot} 
                {}={}
            \frac12
            \begin{bsmallmatrix}
                \I_2 & -\I_2\\
                -\I_2 & \I_2
            \end{bsmallmatrix}
        \).
        Thus, the claim is established by observing that (i) the upper-left $2\times2$ block of $\bar M$ is positive definite, and (ii) the Schur complement of the bottom-right $2\times2$ block of $\bar M$ is given by
        \(
            \frac16
            \begin{bsmallmatrix}
                1 & 2\\
                2 & 4
            \end{bsmallmatrix}
        \),
        which has eigenvalues $0$ and $\nicefrac{5}{6}$.

        \item \ref{it:ex:linear-system:2}:
        By definition, $S$ is $(\mon \proj_{X^\bot}, 0)$\hyp{}semimonotone if and only if
        \[
            \bar{S}
                \coloneqq
            \tfrac12(S + S^\top)
                -
            \mon \proj_{X^\bot}
                =
            \frac12
            \begin{bmatrix}
                -2-\mon & 0 & \mon & 0\\
                0 &-2-\mon & 0 & \mon\\
                \mon & 0 & -\mon & 1\\
                0 & \mon & 1 & -\mon
            \end{bmatrix}
                {}\succeq{}
            0.
        \]
        A necessary condition for $\bar S \succeq 0$ is that its upper-left $2\times2$ block is positive semidefinite.
        Therefore, $S$ is not $(\mon \proj_{X^\bot}, 0)$\hyp{}semimonotone for any $\mon > -2$. 
        Consider the case when $\mon \leq -2$.
        \begin{itemize}
            \item[\(\diamondsuit\)] $\mon = -2$: Then, $\lambda_{\rm min}\bigl(\bar S\bigr) = \frac14(1-\sqrt{17}) < 0$, so $S$ is not $(-2 \proj_{X^\bot}, 0)$\hyp{}semimonotone either.
            \item[\(\diamondsuit\)] $\mon < -2$: Then, the upper-left $2\times2$ block of $\bar S$ is positive definite, and the Schur complement of the bottom-right $2\times2$ block of $\bar S$ is given by
            \(
                \frac12
                \begin{bsmallmatrix}
                    \nicefrac{-2\mon}{(2+\mon)} & 1\\
                    1 & \nicefrac{-2\mon}{(2+\mon)}
                \end{bsmallmatrix}
            \).
            Since the diagonal elements $\nicefrac{-2\mon}{(2+\mon)}$ are strictly negative for any $\mon < -2$, it follows that $S$ is not $(\mon \proj_{X^\bot}, 0)$\hyp{}semimonotone for any $\mon < -2$, completing the proof.
        \end{itemize}
        \item \ref{it:ex:linear-system:3}:
        We begin by showing that \cref{ass:progdec-global} holds.
        \cref{ass:progdec-global:3} is satisfied, as shown in \cref{it:ex:linear-system:1}, since the set
        \(
            \link_X S =
            \Bigl(
                (1,1,1,1),
                (1,-3,-1,3)
            \Bigr)
        \)
        is nonempty and
        \(
            [\mon]_{-}[\com]_{-} = \nicefrac12 < 1
        \).
        \cref{ass:progdec-global:1} also holds since $S$ is continuous and thus outer semicontinuous.
        Finally, the resolvent $J_{\gamma^{-1} S}$ has full domain if and only if the operator $\id + \gamma^{-1} S$ has full range. Given that $S(x) = Mx - m$, this holds if and only if the matrix
        \[
            \I + \gamma^{-1} M
                {}={}
                \begin{bmatrix}
                    \I_2 + \gamma^{-1}M_1 & 0_{2\times 2}\\
                    0_{2\times 2} & \I_2 + \gamma^{-1}M_2
                \end{bmatrix}
        \]
        has full rank.
        Since the determinant of each diagonal block is strictly positive for all $\gamma \in \R$, \cref{ass:progdec-global:2} is established.
        Having shown that \cref{ass:progdec-global} holds, the claim follows directly from \cref{thm:progdec}. \qedhere
    \end{proofitemize}
\end{appendixproof}
\begin{appendixproof}{cor:semimon:lin:functions}
    \begin{proofitemize}
        \item \ref{cor:semimon:lin:functions:global-easy}:
        Let $x, \other{x} \in \pazocal{U}^{\rm link}_1$.
        Owing to the continuity of $F$, it follows from the mean-value theorem \cite[Prop. 7.1.16]{facchinei2003FiniteDimensional} that there exist $n$ points \( \nu_i \) on the line connecting \( x \) and \( \other{x} \) and constants $\eta = (\eta_1, \hdots, \eta_n) \in \Delta_n$ where $\Delta_n$ denotes the $n$-dimensional probability simplex, such that
        \begin{align*}
            F(x) - F(\other{x}) = \sum_{i=1}^n \eta_i \Jac F(\nu_i)(x-\other{x}).
            \numberthis\label{eq:cor:semimon:lin:functions:global:mean-value}
        \end{align*}
        Consequently,
        \ifsiam
\else
            \begin{align*}
                \inner{x - \other{x}, F(x) - F(\other{x})}
                    {}\overrel{\eqref{eq:cor:semimon:lin:functions:global:mean-value}}{}
                \sum_{i=1}^n \eta_i \inner{x - \other{x}, \Jac F(\nu_i)(x-\other{x})}
                    {}\geq{}&
                \mon \|\proj_{X^\bot}(x - \other{x})\|^2
                    +
                \sum_{i=1}^n \eta_i \com \left\|
                    \proj_X \Bigl(F(x) - F(\other{x})\Bigr)
                \right\|^2\\
                    {}={}&
                \mon \|\proj_{X^\bot}(x - \other{x})\|^2
                    +
                \com \left\|
                    \proj_X \Bigl(F(x) - F(\other{x})\Bigr)
                \right\|^2,
            \end{align*}
        \fi
        where the first inequality holds by assumption and the second one holds by Jensen's inequality,
        showing that $F$ is
        \(
            \bigl(
                \mon \proj_{X^\bot},
                \com \proj_X
            \bigr)
        \)\hyp{}semimonotone on $\pazocal{U}^{\rm link}_1 \times \R^n$.
        Finally, the maximality claim for $\com \geq 0$ follows from \Cref{cor:structured-semimon:maximal}.
        
        \item \ref{cor:semimon:lin:functions:global}:
        By \Cref{lem:semimon:lin} the matrices \(\Jac F(\nu_i)\) are
        \(
            \bigl(
                \bar \mon \proj_{X^\bot},
                \com \proj_X
            \bigr)
        \)\hyp{}semimonotone.
        Let $x, \other{x} \in \R^n$.
        Then, it follows from the mean-value theorem that
        \ifsiam
\else
            \begin{align*}
                \inner{x - \other{x}, F(x) - F(\other{x})}
                    {}\overrel{\eqref{eq:cor:semimon:lin:functions:global:mean-value}}{}
                \sum_{i=1}^n \eta_i \inner{x - \other{x}, \Jac F(\nu_i)(x-\other{x})}
                    {}\geq{}&
                \bar \mon \|\proj_{X^\bot}(x - \other{x})\|^2
                    +
                \sum_{i=1}^n \eta_i \com \|\proj_X \Jac F(\nu_i)(x - \other{x})\|^2\\
                    {}\geq{}&
                \bar \mon \|\proj_{X^\bot}(x - \other{x})\|^2
                    +
                \com \left\|
                    \proj_X \left(\sum_{i=1}^n \eta_i \Jac F(\nu_i)(x - \other{x})\right)
                \right\|^2\\
                    {}\overrel{\eqref{eq:cor:semimon:lin:functions:global:mean-value}}{}&
                \bar \mon \|\proj_{X^\bot}(x - \other{x})\|^2
                    +
                \com \left\|
                    \proj_X \Bigl(F(x) - F(\other{x})\Bigr)
                \right\|^2,
            \end{align*}
        \fi
        where the first inequality holds by 
        \(
            \bigl(
                \bar \mon \proj_{X^\bot},
                \com \proj_X
            \bigr)
        \)\hyp{}semimonotonicity of the matrices \(\Jac F(\nu_i)\) and the second one holds by Jensen's inequality,
        showing that $F$ is
        \(
            \bigl(
                \bar \mon \proj_{X^\bot},
                \com \proj_X
            \bigr)
        \)\hyp{}semimonotone, and, since \(\com \geq 0\), maximally so by \Cref{cor:structured-semimon:maximal}.

        \item \ref{cor:semimon:lin:functions:local}: 
        By continuity of $\Jac F$, there exists a convex neighborhood $\pazocal{U}^{\rm link}_1$ of $\tilde{x}$ and a function $\varepsilon : \pazocal{U}^{\rm link}_1 \rightarrow [0, \mon)$
        satisfying
        $\varepsilon(\tilde{x}) = 0$
        and
        \(
            \sup_{y \in \pazocal{U}^{\rm link}_1} \varepsilon(y) < \mon
        \)
        such that
        \(
            \|\proj_X \bigl(\Jac F(\tilde{x}) - \Jac F(y)\bigr)\proj_X\| \leq \varepsilon(y)
        \)
        for all \(y \in \pazocal{U}^{\rm link}_1\).
        Let $x \in X$ and $y \in \pazocal{U}^{\rm link}_1$. Then, it holds that
        \begin{align*}
            \langle x, \Jac F(y) x \rangle
                {}={}
            \langle x, \Jac F(\tilde{x}) x + (\Jac F(y) - \Jac F(\tilde{x}) ) x\rangle
                {}\geq{}&
            \langle x, \Jac F(\tilde{x}) x \rangle - \varepsilon(y) \|x\|^2\\
                {}\overrel[\geq]{\eqref{cor:semimon:lin:functions:local:ass}}{}&
            \left(\mon - \varepsilon(y)\right) \|x\|^2 + \com \|\proj_X \Jac F(\tilde{x}) x\|^2.
            \numberthis\label{cor:semimon:lin:functions:local:res1}
        \end{align*}
        For any function      
        \(
            \nu : \pazocal{U}^{\rm link}_1 \rightarrow (0, 1]
        \),
        Young's inequality gives
        \begin{align*}
            \|\proj_X \Jac F(\tilde{x}) x\|^2
                {}={}&
            \left\|\proj_X \Jac F(y) x\right\|^2
                {}+{}
            \ifsiam
\else
                2\left\langle \proj_X \Jac F(y) x, \proj_X \Bigl(\Jac F(\tilde{x}) - \Jac F(y) \Bigr)x 
                \right\rangle
            \fi
                {}+{}
            \left\|\proj_X \Bigl(\Jac F(\tilde{x}) - \Jac F(y) \Bigr)x\right\|^2\\
                {}\geq{}&
            \left(1-\nu(y)\right)\left\|\proj_X \Jac F(y) x\right\|^2
                {}+{}
            (1-\tfrac1{\nu(y)})\left\|\proj_X \Bigl(\Jac F(\tilde{x}) - \Jac F(y) \Bigr)x\right\|^2
            \\
                {}\geq{}&
            \left(1-\nu(y)\right)\left\|\proj_X \Jac F(y) x\right\|^2
                {}+{}
            (1-\tfrac1{\nu(y)})\varepsilon(y)^2\left\|x\right\|^2,
        \end{align*}
        Combining this with \eqref{cor:semimon:lin:functions:local:res1} and using that $\com \geq 0$, this implies that
        \begin{align*}
            \langle x, \Jac F(y) x \rangle
                {}\geq{}
            \underbrace{\bigl(\mon - \varepsilon(y) + (1-\tfrac1{\nu(y)})\varepsilon(y)^2\com\bigr)}_{\eqqcolon\, \mon(y)} \|x\|^2 +
            \underbrace{\left(1-\nu(y)\right)\com}_{\eqqcolon\, \com(y)} \left\|\proj_X \Jac F(y) x\right\|^2.
            \numberthis\label{eq:cor:semimon:lin:functions:local:intermediate}
        \end{align*}
        Note that by construction
        \(
            \bar \com
                {}\coloneqq{}
            \inf_{y \in \pazocal{U}^{\rm link}_1} \com(y)
                {}\geq{}
            0.
        \)
        Now, without loss of generality, let the function $\nu$ satisfy
        \[
            \inf_{y \in \pazocal{U}^{\rm link}_1} \nu(y)
                >
            \inf_{y \in \pazocal{U}^{\rm link}_1}
            \frac{\varepsilon(y)^2\com}{\mon - \varepsilon(y) + \varepsilon(y)^2\com}
                \in
            (0,1),
        \]
        ensuring that $\inf_{y \in \pazocal{U}^{\rm link}_1} \mon(y) > 0$.
        Then, the claim follows from \Cref{cor:semimon:lin:functions:global}.\qedhere
    \end{proofitemize}
\end{appendixproof}

\end{appendix}

\ifsiam
\else
  \bibliographystyle{arxivplain}%
  \bibliography{TeX/Minty-DRS.bib}%
\fi

\end{document}